\documentclass[11pt]{article}
\usepackage{amssymb}
\usepackage{amsmath}
\usepackage{colortbl,xspace}

\newtheorem{theo}{Theorem}
\newtheorem{lem}{Lemma}
\newtheorem{prop}{Proposition}
\newtheorem{cor}{Corollary}
\newtheorem{rem}{Remark}
\newtheorem{exa}{Example}
\newtheorem{definiti}{Definition}

\newenvironment{dem}[1][Proof]{\noindent \textbf{#1.} }{\ \rule{0.5em}{0.5em}}
\input{tcilatex}

\begin{document}

\title{Lispchitz modulus of the argmin mapping in convex quadratic
optimization\thanks{%
This research has been partially supported by Grant PID2022-136399NB-C22
from MICINN, Spain, and ERDF, "A way to make Europe", European Union. This
research was also partially supported by the CIPROM/2024/34 grant, funded by
the Conselleria de Educaci\'{o}n, Cultura, Universidades y Empleo,
Generalitat Valenciana.}}
\date{}
\author{M.J. C\'{a}novas\thanks{%
Center of Operations Research, Miguel Hern\'{a}ndez University of Elche,
03202 Elche (Alicante), Spain;\ canovas@umh.es, parra@umh.es.} \and M.
Fukushima\thanks{%
The Kyoto College of Graduate Studies for Informatics, Kyoto, Japan;
fukushima\_masao@yahoo.co.jp.} \and J. Parra\footnotemark[2] \and {\small {%
Dedicated to Prof. Marco A.L\'{o}pez as a token}} \and {\small {of
recognition for his outstanding career}}}
\maketitle

\begin{abstract}
This paper was initially motivated by the computation of the Lipschitz
modulus of the metric projection on polyhedral convex sets in the Euclidean
space when both the reference point and the polyhedron where it is projected
are subject to perturbations. The paper tackles the more general problem of
computing the Lipschitz modulus of the argmin mapping in the framework of
canonically perturbed convex quadratic problems. We point out the fact that
a point-based formula (depending only on the nominal data) for such a
modulus is provided. In this way, the paper extends to the current quadratic
setting some results previously developed in linear programming. As an
application, we provide a point-based formula for the Lipschitz modulus of
the metric projection on a polyhedral convex set.


\bigskip

\end{abstract}

\section{Introduction}

The initial motivation of this paper was to compute the \emph{Lipschitz
modulus }of the \emph{metric projection function }on varying polyhedral
convex sets in the $n$-dimensional Euclidean space. Specifically, we
consider a parametrized family of polyhedral convex sets $\left\{ \mathcal{F}%
\left( b\right) :b\in \mathbb{R}^{m}\right\} $ given by%
\begin{equation}
\mathcal{F}\left( b\right) :=\left\{ x\in \mathbb{R}^{n}:a_{i}^{\prime
}x\leq b_{i},\text{ }i\in I:=\left\{ 1,...,m\right\} \right\} ,\text{ }%
b=\left( b_{i}\right) _{i\in I}\in \mathbb{R}^{m},  \label{eq_F}
\end{equation}%
where $a_{i}\in \mathbb{R}^{n},$ $i=1,...,m,$ are fixed. Elements of $%
\mathbb{R}^{k}$,$\ $for any $k\in \mathbb{N},$ are considered as column
vectors and the prime denotes transposition, so that $a_{i}^{\prime }x$ is
the usual inner product of $a_{i}$ and $x$ in $\mathbb{R}^{n}$. The metric
projection function $\mathcal{P}:\mathbb{R}^{n}\times \mathbb{R}%
^{m}\rightarrow \mathbb{R}^{n}$ is given by%
\begin{equation}
\mathcal{P}\left( z,b\right) :=\arg \min \left\{ \left\Vert z-x\right\Vert
_{2}:x\in \mathcal{F}\left( b\right) \right\} ,  \label{eq_P}
\end{equation}%
where $\left\Vert \cdot \right\Vert _{2}$ represents the Euclidean norm in $%
\mathbb{R}^{n}.$ Given a fixed $\overline{z}\in \mathbb{R}^{n},$ it is
well-known that $\mathcal{P}\left( \overline{z},b\right) $ consists of the
unique optimal solution of the strictly convex quadratic optimization
problem given by minimizing $x^{\prime }x-2\overline{z}^{\prime }x$ subject
to $x\in \mathcal{F}\left( b\right) $.

One can find in the literature different contributions on the Lipschitz
behavior of the metric projection in different frameworks of perturbations
(some of them in infinite-dimensional spaces);\ see, e.g. \cite{Ab79,
BaGo12, BeRu19, Yen95}. In particular, Yen \cite{Yen95} establishes the
Lipschitz property of function $\mathcal{P},$ defined in (\ref{eq_P}), at a
nominal (fixed) element $\left( \overline{z},\overline{b}\right) \in \mathbb{%
R}^{n}\times \mathbb{R}^{m},$ provided that $\mathcal{F}\left( b\right) \neq
\emptyset $ for $b\in \mathbb{R}^{m}$ close enough to $\overline{b}.$ That
paper ensures the existence of a certain \emph{Lipschitz constant} for $%
\mathcal{P}$ at $\left( \overline{z},\overline{b}\right) .$ The present
paper provides, as a byproduct, an implementable formula for computing the
infimum of all Lipschitz constants (the so-called \emph{Lipschitz modulus})%
\emph{\ }of $\mathcal{P}$ at $\left( \overline{z},\overline{b}\right) $; see
Section 7.

In this paper we deal with a more general framework of canonically perturbed
convex quadratic optimization problems, which includes ordinary LP problems.
Formally, we consider the \emph{optimal set mapping, }$\mathcal{S}:\mathbb{R}%
^{n}\times \mathbb{R}^{m}\rightrightarrows \mathbb{R}^{n},$\emph{\ }defined
by 
\begin{equation}
\mathcal{S}\left( c,b\right) :=\arg \min \left\{ \frac{1}{2}x^{\prime
}Qx+c^{\prime }x:x\in \mathcal{F}\left( b\right) \right\} ,  \label{eq_S}
\end{equation}%
where $Q$ is a symmetric positive semidefinite $n\times n$ matrix, and $%
\mathcal{F}:\mathbb{R}^{m}\rightrightarrows \mathbb{R}^{n}$ is the \emph{%
feasible set mapping }defined in (\ref{eq_F}). For convenience, let us
denote by $A$ the (fixed) matrix whose $i$-th row is $a_{i}^{\prime },$ $%
i\in I.$ In this framework $\left( c,b\right) \in \mathbb{R}^{n}\times 
\mathbb{R}^{m}$ is the parameter to be perturbed. Thus, we are placed in the
setting of canonical perturbations: linear perturbations of the objective
function together with perturbations of the right-hand side (RHS) of the
constraints. For the technical developments of this paper we require $Q$ to
be fixed. The continuity properties of $\mathcal{S}$ in this context have
been analyzed in \cite[Section 5.5]{BGKKT82}. A model where matrix $Q$ can
be perturbed is considered in \cite[Section 5.3]{BGKKT82}; see also \cite%
{LeeYen14}, where the feasible set is the closed Euclidean ball centered at
the origin and whose radius may be perturbed.

The main aim of this work consists of analyzing the Lipschitzian behavior of
the \emph{argmin mapping}, $\mathcal{S}$, and providing point-based
expressions for the modulus (i.e., given in terms of the nominal data, not
involving data in a neighborhood). Specifically, we analyze the \emph{Aubin
property} (also called pseudo-Lipschitz or Lipschitz-Like) of $\mathcal{S}$
at $\left( \left( \overline{c},\overline{b}\right) ,\overline{x}\right) \in 
\mathrm{gph}\mathcal{S}$ (the graph of $\mathcal{S)}.$ Roughly speaking, the
fulfilment of this property ensures the existence of upper bounds (Lipschitz
constants) on the rate of variation of optimal solutions around $\overline{x}
$ with respect to data perturbations (i.e., perturbations of $\left( 
\overline{c},\overline{b}\right) )$; see Section 2.1 for the formal
definitions.

The Aubin property of a multifunction is a well-known notion of variational
analysis, with many applications to optimization (to both, theory and
algorithms). The reader is addressed to the monographs \cite{DoRo, Ioffe17,
KlKu02, mor06a, rw} for a comprehensive analysis of this and other
Lipschitz-type and regularity properties. For instance, some applications of 
\emph{calmness} and \emph{error bounds} to the speed of convergence of
certain algorithms are given in \cite{KYF04, KK09}. At this moment we
underline the role played by the \emph{strong Lipschitz stability }and a 
\emph{certain Kojima-type }stability condition (see, e.g., \cite[Chapter 8]%
{KlKu02} and \cite{KlaKu05}). The existence of Lipschitz-type constants for
more general quadratic models, not requiring $Q$ in (\ref{eq_S}) to be
positive semidefinite and with a polyhedral feasible set mapping (its graph
is a finite union of convex polyhedra), has been studied in the recent paper 
\cite{Klatte25}.

The direct antecedent of this paper can be found in the context of linear
optimization problems, where $Q$ is the null matrix, i.e., $Q=0_{n\times n}$%
. In the linear framework the Lipschitz modulus of $\mathcal{S}$ at $\left(
\left( \overline{c},\overline{b}\right) ,\overline{x}\right) $ has already
been determined through exact formulae in \cite{CGP08}. The current setting
of convex quadratic problems is notably different from the linear one. To
start with, when confined to the linear setting, the Aubin property is
equivalent to the so called \emph{N\"{u}rnberger Condition} (NC for short),
which, roughly speaking, forces the Karush-Kuhn-Tucker (KKT) conditions to
hold at $\overline{x}\in \mathcal{S}\left( \overline{c},\overline{b}\right) $
with at least $n$ active constraints (see Section 2.2 for details). However,
the NC is shown to be very restrictive in the convex quadratic context (see
Section 2.2). In \cite{CKLP07}, the extension of NC for convex semi-infinite
problems, called in later papers \emph{Extended N\"{u}rnberger Condition }%
(ENC)\emph{,} was shown to be a sufficient condition for the Aubin property
of the corresponding mapping $\mathcal{S}$. After that work, several steps
approaching the Lipschitz modulus in the convex setting were given under the
ENC (see, for instance, \cite{CHLP08}). In contrast, the current work
tackles the Aubin property of $\mathcal{S}$ without requiring the NC. As a
counterpart, the family of \emph{minimal KKT} \emph{sets of active indices }%
becomes a key tool (see Section 2.2 for the definitions and first results).
We point out that this family of sets of indices was introduced in \cite%
{chlp16}, playing a key role in the computation of the calmness modulus of
the optimal set (argmin) mapping for linear programs.

The structure of the paper is as follows: After the section of preliminaries
and first results, Section 3 introduces the minimal KKT subproblems
associated with the minimal KKT sets of active indices. The Aubin property
of such subproblems is analyzed in Section 4. Proposition \ref{Prop>=} in
that section and the subsequent Example \ref{Exa1} show that the minimal
sets are not enough to determine the Lipschitz modulus of $\mathcal{S}.$ An
extended class of index subsets is introduced and analyzed in Section 5.
With this tool, Section 6 characterizes the Aubin property of $\mathcal{S}$
and provides the aimed point-based expression for its Lipschitz modulus. We
finish the paper with a short section of conclusions, when we revisit the
problem of the metric projection

\section{Preliminaries and first results}

Throughout the paper, the space of variables, $\mathbb{R}^{n},$ is equipped
with an arbitrary norm, $\left\Vert \cdot \right\Vert ,$ whose corresponding 
\emph{dual norm} is given by $\left\Vert x\right\Vert _{\ast
}=\max_{\left\Vert y\right\Vert \leq 1}\left\vert x^{\prime }y\right\vert ,$
while the parameter space $\mathbb{R}^{n}\times \mathbb{R}^{m}$ is endowed
with the norm 
\begin{equation}
\left\Vert \left( c,b\right) \right\Vert :=\max \left\{ \left\Vert
c\right\Vert _{\ast },\left\Vert b\right\Vert _{\infty }\right\} ,\text{
where }\left\Vert b\right\Vert _{\infty }=\max_{1\leq i\leq m}\left\vert
b_{i}\right\vert .  \label{eq_norms}
\end{equation}%
For the sake of simplicity in the notation, distances in both the space of
variables $\mathbf{(}\mathbb{R}^{n})$ and the parameter space $\mathbf{(}%
\mathbb{R}^{n}\times \mathbb{R}^{m})$ will be denoted by the same symbol $d,$
and they will be distinguished by the context. As usual, the distance from a
point $z\in \mathbb{R}^{n}$\ to a subset $S\subset \mathbb{R}^{n}$ is given
by%
\begin{equation*}
d\left( z,S\right) :=\inf \{d\left( z,x\right) :x\in S\},
\end{equation*}%
with the convention $d\left( x,\emptyset \right) :=+\infty $ (according to $%
\inf \emptyset :=+\infty $).

Given $X\subset \mathbb{R}^{k},$ $k\in \mathbb{N},$ we denote by \textrm{conv%
}$X$ and \textrm{cone}$X$ the \emph{convex hull} and the \emph{conical
convex hull }of $X$, respectively. It is assumed that \textrm{cone}$X$
always contains the zero-vector $0_{k}$, in particular \textrm{cone}$%
(\emptyset )=\{0_{k}\}.$ If $X$ is a subset of any topological space, $%
\mathrm{int}X,$ $\mathrm{cl}X$ and $\mathrm{bd}X$ stand, respectively, for
the interior, the closure, and the boundary of $X.$ As usual, $\mathbb{R}%
_{+}=[0,+\infty \lbrack .$

\subsection{On the Lipschitz behavior of the optimal set}

Here we recall the Aubin property specified directly for our optimal set
mapping $\mathcal{S}$ defined in (\ref{eq_S}), although the definition is
the same for any set-valued mapping between metric spaces. $\mathcal{S}$ has
the \emph{Aubin property }at $\left( \left( \overline{c},\overline{b}\right)
,\overline{x}\right) \in \mathrm{gph}\mathcal{S}$ if there exist
neighborhoods $U\subset \mathbb{R}^{n}$ and $V\subset \mathbb{R}^{n}\times 
\mathbb{R}^{m},$ of $\overline{x}$ and $\left( \overline{c},\overline{b}%
\right) ,$ respectively, and a constant $\kappa \geq 0$ such that 
\begin{equation}
d\left( x^{2},\mathcal{S}\left( c^{1},b^{1}\right) \right) \leq \kappa
d\left( \left( c^{2},b^{2}\right) ,\left( c^{1},b^{1}\right) \right) ,\text{ 
}  \label{eq_def_Lips}
\end{equation}%
whenever $x^{2}\in \mathcal{S}\left( c^{2},b^{2}\right) \cap U$ and $\left(
c^{1},b^{1}\right) ,\left( c^{2},b^{2}\right) \in V.$ The infimum of all
(Lipschitz) constants $\kappa $ such that (\ref{eq_def_Lips}) holds for some
associated neighborhoods $U$ and $V$ is the Lipschitz modulus of $\mathcal{S}
$ at $\left( \left( \overline{c},\overline{b}\right) ,\overline{x}\right) ,$ 
$\mathrm{lip}\mathcal{S}\left( \left( \overline{c},\overline{b}\right) ,%
\overline{x}\right) .$ Alternatively, we can write 
\begin{equation*}
\mathrm{lip}\mathcal{S}\left( \left( \overline{c},\overline{b}\right) ,%
\overline{x}\right) =\underset{x^{2}\rightarrow \overline{x},\text{ }%
x^{2}\in \mathcal{S}\left( c^{2},b^{2}\right) }{\underset{\left(
c^{1},b^{1}\right) ,\left( c^{2},b^{2}\right) \rightarrow \left( \overline{c}%
,\overline{b}\right) ,}{\lim \sup }}\frac{d\left( x^{2},\mathcal{S}\left(
c^{1},b^{1}\right) \right) }{d\left( \left( c^{2},b^{2}\right) ,\left(
c^{1},b^{1}\right) \right) },
\end{equation*}%
under the convention $\frac{0}{0}:=0.$ We write $\mathrm{lip}\mathcal{S}%
\left( \left( \overline{c},\overline{b}\right) ,\overline{x}\right) =+\infty 
$ when $\mathcal{S}$ has not the Aubin property at $\left( \left( \overline{c%
},\overline{b}\right) ,\overline{x}\right) $. The \emph{calmness property }\
of $\mathcal{S}$ at $\left( \left( \overline{c},\overline{b}\right) ,%
\overline{x}\right) $ and the corresponding calmness modulus, $\mathrm{clm}%
\mathcal{S}\left( \left( \overline{c},\overline{b}\right) ,\overline{x}%
\right) ,$ come from fixing $\left( c^{1},b^{1}\right) =\left( \overline{c},%
\overline{b}\right) $ in the previous definitions.

$\mathcal{S}$ is \emph{strongly Lipschitz stable }at $\left( \overline{c},%
\overline{b}\right) $ if $\mathcal{S}$ is single-valued in some neighborhood
of $\left( \overline{c},\overline{b}\right) $ and Lipschitz continuous at $%
\left( \overline{c},\overline{b}\right) .$ The single-valuedness of $%
\mathcal{S}$ around $\left( \overline{c},\overline{b}\right) $ together with
its continuity at $\left( \overline{c},\overline{b}\right) $ is usually
referred to as the \emph{Kojima-type stability}.

In the following result, the equivalence between $\left( i\right) $ and $%
\left( ii\right) $ holds, in fact, in a more general framework of
optimization problems in Hilbert spaces with linearly perturbed objective
functions and arbitrary feasible sets (see \cite[Corollary 4.7]{KlKu02}).
The remaining equivalences were established in \cite[Theorem 16]{CKLP07}.

\begin{prop}
\label{Prop1}Let $\left( \left( \overline{c},\overline{b}\right) ,\overline{x%
}\right) \in \mathrm{gph}\mathcal{S}.$ The following conditions are
equivalent:

$\left( i\right) $ $\mathcal{S}$ has the Aubin property at $\left( \left( 
\overline{c},\overline{b}\right) ,\overline{x}\right) ;$

$\left( ii\right) $ $\mathcal{S}$ is strongly Lipschitz stable\emph{\ }at $%
\left( \overline{c},\overline{b}\right) ;$

In the linear case, when $Q=0_{n\times n}$, we can add the following:

$\left( iii\right) $ $\mathcal{S}$ is single-valued and continuous in some
neighborhood of $\left( \overline{c},\overline{b}\right) ;$

$\left( iv\right) $ $\mathcal{S}$ is single-valued in some neighborhood of $%
\left( \overline{c},\overline{b}\right) .$
\end{prop}

The next proposition is focused on the Aubin property of the feasible set
mapping $\mathcal{F}.$ The computation of $\mathrm{lip}\mathcal{F}\left( 
\overline{b},\overline{x}\right) $ was addressed in \cite[Corollary 3.2]%
{CDLP05}. From now on, $\mathrm{dom}\mathcal{F}$ represents the domain of $%
\mathcal{F};$ i.e., 
\begin{equation*}
\mathrm{dom}\mathcal{F}:=\left\{ b\in \mathbb{R}^{m}:\mathcal{F}\left(
b\right) \neq \emptyset \right\} .
\end{equation*}%
Recall that Slater constraint qualification (SCQ) holds at $\overline{b}$ if
there exists $y\in \mathbb{R}^{n}$ satisfying $a_{i}^{\prime }y<\overline{b}%
_{i},$ $i=1,...,m;$ such a point $y$ is referred to as a Slater point of the
linear system associated with $\overline{b}.$ Moreover, given $\left( 
\overline{b},\overline{x}\right) \in \mathrm{gph}\mathcal{F},$ $I_{\overline{%
b}}\left( \overline{x}\right) $ denotes the set of active indices at $%
\overline{x};$ i.e., 
\begin{equation*}
I_{\overline{b}}\left( \overline{x}\right) :=\{i\in
\{1,...,m\}:a_{i}^{\prime }\overline{x}=\overline{b}_{i}\}.
\end{equation*}%
In the next proposition, `$\left( i\right) \Leftrightarrow $ $\left(
ii\right) $' comes from the classical Robinson-Ursescu theorem (for
multifunctions with closed and convex graphs), while equivalences `$\left(
ii\right) \Leftrightarrow \left( iii\right) \Leftrightarrow \left( iv\right) 
$' can be found, for instance, in \cite[Theorem 6.1]{libro}.

\begin{prop}
\label{Prop_SCQ}Let $\left( \overline{b},\overline{x}\right) \in \mathrm{gph}%
\mathcal{F}.$ The following conditions are equivalent:

$\left( i\right) $ $\mathcal{F}$ has the Aubin property at $\left( \overline{%
b},\overline{x}\right) ;$

$\left( ii\right) $ $\overline{b}\in \mathrm{intdom}\mathcal{F}$ $;$

$\left( iii\right) $ SCQ holds at $\overline{b};$

$\left( iv\right) $ $0_{n}\notin \mathrm{conv}\{a_{i}:i\in I_{\overline{b}%
}\left( \overline{x}\right) \}.$
\end{prop}

Clearly condition $\left( ii\right) $ above is necessary for the Aubin
property of $\mathcal{S}$ at $\left( \left( \overline{c},\overline{b}\right)
,\overline{x}\right) \in \mathrm{gph}\mathcal{S}$.

\textbf{Assumption: }Throughout the paper we assume that the equivalent
conditions of Proposition \ref{Prop_SCQ} hold at our nominal $\left( \left( 
\overline{c},\overline{b}\right) ,\overline{x}\right) \in \mathrm{gph}%
\mathcal{S}.$

\begin{rem}
\label{Rem FW}\emph{Under the assumption above, the only possibility to have 
}$\mathcal{S}\left( c,b\right) =\emptyset $ \emph{for }$\left( c,b\right) $ 
\emph{close enough to} $\left( \overline{c},\overline{b}\right) $\emph{\ is
that the objective function} $x\mapsto \frac{1}{2}x^{\prime }Qx+c^{\prime }x$
\emph{is unbounded from below on }$\mathcal{F}\left( b\right) .$ \emph{This
fact comes, for instance, from the classical Frank-Wolfe Theorem.}
\end{rem}

\subsection{N\"{u}rnberger Condition \ and minimal KKT sets of indices}

Recall that $\left( \left( c,b\right) ,x\right) \in \mathrm{gph}\mathcal{S}$
if and only if $x$ is a KKT point for the convex quadratic optimization
problem associated with $\left( c,b\right) $ (recall that $Q$ is positive
semidefinite), which can be expressed as the existence of scalars $\lambda
_{i}\geq 0,$ $i=1,...,m,$ such that 
\begin{equation}
\left\{ 
\begin{array}{c}
-\left( Qx+c\right) =\dsum\limits_{i=1}^{m}\lambda _{i}a_{i}, \\ 
\lambda _{i}\left( a_{i}^{\prime }x-b_{i}\right) =0,i=1,...,m, \\ 
a_{i}^{\prime }x-b_{i}\leq 0,i=1,...,m.%
\end{array}%
\right.  \label{eq_kkt}
\end{equation}%
In this way $\left( \left( c,b\right) ,x,\lambda \right) ,$ with $\lambda
=\left( \lambda _{i}\right) _{i=1,...,m},$ can be seen as the solution of a
linear system with complementarity constraints. We express the KKT
conditions (\ref{eq_kkt}) by using the set of active indices: $\left( \left(
c,b\right) ,x\right) \in \mathrm{gph}\mathcal{S}$ if and only if $x\in 
\mathcal{F}\left( b\right) $ and 
\begin{equation*}
-(Qx+c)\in \mathrm{cone}\left\{ a_{i},i\in I_{b}\left( x\right) \right\} .
\end{equation*}%
Adapting the definition of the N\"{u}rnberger Condition (NC) (given in \cite[%
Section 4]{CKLP07} for linear problems) to our current quadratic setting
(see also condition (10) in \cite[Section 4]{CKLP07} for the ENC referred to
in Section 1), we say that $\left( \left( c,b\right) ,x\right) \in \mathrm{%
gph}\mathcal{S}$ satisfies the NC if $b$ satisfies the SCQ and 
\begin{equation*}
\left\{ 
\begin{array}{c}
D\subset I_{b}\left( x\right) \\ 
-(Qx+c)\in \mathrm{cone}\left\{ a_{i},i\in D\right\}%
\end{array}%
\right\} \Rightarrow \left\vert D\right\vert \geq n,
\end{equation*}%
where $\left\vert D\right\vert $ stands for the cardinality of $D.$

The following result which comes from \cite[Theorems 10 and 16]{CKLP07}
establishes that the NC at $\left( \left( c,b\right) ,x\right) $ is a
sufficient condition for the Aubin property of $\mathcal{S}$ at $\left(
\left( \overline{c},\overline{b}\right) ,\overline{x}\right) .$ The same
paper points out that the NC is not necessary for the Aubin property of $%
\mathcal{S}$ (\cite[Remark 11]{CKLP07}).

\begin{prop}
\label{Prop NC linear}Let $\left( \left( c,b\right) ,x\right) \in \mathrm{gph%
}\mathcal{S}.$ If $\left( \left( c,b\right) ,x\right) $ satisfies the NC,
then $\mathcal{S}$ has the Aubin property at $\left( \left( c,b\right)
,x\right) .$ The converse implication holds true in the case when $%
Q=0_{n\times n}.$
\end{prop}

Now we introduce the (already mentioned) key tool used in this work
(specifically in Theorem \ref{Th main}) for characterizing the Aubin
property of $\mathcal{S},$ namely, the minimal KKT set of active indices. As
advanced in Section 1, this tool was already used in \cite{chlp16} for
computing the calmness modulus of $\mathcal{S}$, provided that $Q=0_{n\times
n}.$ Here, we extend the definition to convex quadratic problems. As also
mentioned at the end of Section 1, this tool is not enough to determine the
Lipschitz modulus of $\mathcal{S}.$

\begin{definiti}
Let $\left( \left( c,b\right) ,x\right) \in \mathrm{gph}\mathcal{S}.$ We say
that $D\subset I_{b}\left( x\right) $ is a \emph{minimal KKT set of indices}
at $\left( \left( c,b\right) ,x\right) $ if 
\begin{equation}
-(Qx+c)\in \mathrm{cone}\left\{ a_{i},i\in D\right\}  \label{eq_D}
\end{equation}%
and $D$ is minimal with respect to the inclusion order (i.e., there is no $%
\widetilde{D}\subsetneqq D$ satisfying (\ref{eq_D})). The family of all 
\emph{minimal KKT sets of indices} at $\left( \left( c,b\right) ,x\right) $
is denoted by $\mathcal{M}_{c,b}\left( x\right) .$
\end{definiti}

\begin{rem}
Observe that $\emptyset \in \mathcal{M}_{c,b}\left( x\right) ,$ with $\left(
\left( c,b\right) ,x\right) \in \mathrm{gph}\mathcal{S},$ means $Qx+c=0_{n}.$
\end{rem}

\begin{lem}
\label{Lem1}Let $\left( \left( c,b\right) ,x\right) \in \mathrm{gph}\mathcal{%
S}.$ Then one has:

$\left( i\right) $ If $D\in \mathcal{M}_{c,b}\left( x\right) ,$ then the set
of vectors $\left\{ a_{i},i\in D\right\} \subset \mathbb{R}^{n}$ is linearly
independent;

$\left( ii\right) $ The NC holds at $\left( \left( c,b\right) ,x\right) $ if
and only if 
\begin{equation}
\mathcal{M}_{c,b}\left( x\right) =\left\{ D\subset I_{b}\left( x\right)
:-(Qx+c)\in \mathrm{cone}\left\{ a_{i},i\in D\right\} ,\text{ }\left\vert
D\right\vert =n\right\} .  \label{eq_M_UnderNC}
\end{equation}%
Condition `$\left\vert D\right\vert =n$' in (\ref{eq_M_UnderNC}) can be
replaced with `$\left\{ a_{i},i\in D\right\} $ is a basis of $\mathbb{R}^{n}$%
'.
\end{lem}

\begin{dem}
$\left( i\right) $ Consider the nontrivial case $D\neq \emptyset .$
Reasoning by contradiction, if $\left\{ a_{i},i\in D\right\} $ is linearly
dependent, we follow the standard argument of Carath\'{e}odory's theorem to
find $\widetilde{D}\subsetneqq D$ such that $-(Qx+c)\in \mathrm{cone}\left\{
a_{i},i\in \widetilde{D}\right\} ,$ which contradicts the minimality of $D.$

$\left( ii\right) $ Clearly (\ref{eq_M_UnderNC}) implies the NC. Under the
NC, $D\in \mathcal{M}_{c,b}\left( x\right) $ implies $\left\vert
D\right\vert \geq n,$ which together with condition $\left( i\right) ,$
yields that $\left\{ a_{i},i\in D\right\} $ is a basis of $\mathbb{R}^{n}.$

To establish the converse, take $D\subset I_{b}\left( x\right) ,$ such that $%
-(Qx+c)\in \mathrm{cone}\left\{ a_{i},i\in D\right\} $ and $\left\vert
D\right\vert =n.\,$If\ $D\notin \mathcal{M}_{c,b}\left( x\right) ,$ there
would exist $\widetilde{D}\subsetneqq D$ such that $-(Qx+c)\in \mathrm{cone}%
\left\{ a_{i},i\in \widetilde{D}\right\} .$ Since $\left\vert \widetilde{D}%
\right\vert <\left\vert D\right\vert =n,$ we obtain a contradiction with the
NC.
\end{dem}

From now on, given any $D\subset I$ ($=\{1,...,m\}$), we denote by $A_{D}$
the submatrix of $A$ formed by rows $a_{i}^{\prime },$ with $i\in D$ (just
deleting the remaining ones), and for $b\in \mathbb{R}^{m}$ we \ write $%
b_{D}=\left( b_{i}\right) _{i\in D}\in \mathbb{R}^{D}$. For convenience in
the notation, when we delete some constraints we will not relabel the
remaining ones, so that we write $\mathbb{R}^{D}$ instead of $\mathbb{R}%
^{\left\vert D\right\vert }$.

In the linear case $Q=0_{n\times n},$ the family of subsets in the right
member of (\ref{eq_M_UnderNC}) was already used in \cite{CGP08} as a key
tool for deriving an implementable formula for $\mathrm{lip}\mathcal{S}%
\left( \left( \overline{c},\overline{b}\right) ,\overline{x}\right) .$ With
the current notation, condition $\left( ii\right) $ in the previous lemma
allows us to write that formula as follows.

\begin{theo}[{\protect\cite[Corollary 2]{CGP08}}]
\label{Th_CGP08}Assume that $Q=0_{n\times n}$ and that $\left( \left( 
\overline{c},\overline{b}\right) ,\overline{x}\right) \allowbreak \in 
\mathrm{gph}\mathcal{S}$ satisfies the NC. Then,%
\begin{equation*}
\mathrm{lip}\mathcal{S}\left( \left( \overline{c},\overline{b}\right) ,%
\overline{x}\right) =\max_{D\in \mathcal{M}_{\bar{c},\bar{b}}\left( \bar{x}%
\right) }\left\Vert A_{D}^{-1}\right\Vert .
\end{equation*}
\end{theo}

Just for completeness we recall some alternative expressions for $\left\Vert
A_{D}^{-1}\right\Vert $\ traced out from \cite[Section 4]{CGP08}: 
\begin{align*}
\left\Vert A_{D}^{-1}\right\Vert & :=\max\limits_{\beta \in \mathbb{R}%
^{D},~\left\Vert \beta \right\Vert _{\infty }\leq 1}\left\Vert
A_{D}^{-1}\beta \right\Vert =\left( \min_{\left\Vert \lambda \right\Vert
_{1}=1}\left\Vert A_{D}^{\prime }\lambda \right\Vert \right) ^{-1} \\
& =d_{\ast }\left( 0_{n},\limfunc{bd}\limfunc{conv}\left\{ \pm a_{i},i\in
D\right\} \right) ^{-1},
\end{align*}%
provided that $A_{D}$ is\ a\ square\ and\ invertible\ matrix, $A_{D}^{\prime
}$ denotes the transpose of $A_{D},$ and $\left\Vert \cdot \right\Vert _{1}$
stands for the $l^{1}$ norm in $\mathbb{R}^{D}.$

\section{Minimal KKT subproblems}

For convenience we translate the KKT conditions and related concepts to the
matrix form. With this notation, the KKT conditions (\ref{eq_kkt}) can be
written as follows: $\left( \left( c,b\right) ,x\right) \in \mathrm{gph}%
\mathcal{S}$ if and only if there exist $D\subset \{1,...,m\}$ and $\lambda
\in \mathbb{R}^{D}$ such that 
\begin{equation}
\left\{ 
\begin{array}{l}
\left( 
\begin{array}{cc}
Q & A_{D}^{\prime } \\ 
A_{D} & 0_{D\times D}%
\end{array}%
\right) \dbinom{x}{\lambda }=\dbinom{-c}{b_{D}}, \\ 
Ax\leq b,\text{ }\lambda \geq 0_{D}.%
\end{array}%
\right.  \label{eq_kkt2}
\end{equation}

\begin{rem}
\label{Remark_empty}$\left( i\right) $ \emph{For the sake of simplicity (in
order to avoid distinguishing cases), we do not exclude the possibility }$%
D=\emptyset $ \emph{in (\ref{eq_kkt2}).\thinspace\ In such a case, }$%
A_{\emptyset },$ $A_{\emptyset }^{\prime }$ \emph{and }$0_{\emptyset \times
\emptyset }\ $\emph{will be understood as `empty matrices', and }$%
b_{\emptyset }$ \emph{and }$\lambda $\emph{\ as `empty vectors', so that
they do not add any row and column to the corresponding partitioned
matrices. In this way, we have that} 
\begin{equation}
\left( \left( c,b\right) ,x\right) \in \mathrm{gph}\mathcal{S\Leftrightarrow 
}(\text{\ref{eq_kkt2}})\emph{\,}\text{\emph{is satisfied for some }}D\subset
\{1,...,m\}\text{ \emph{and }}\lambda \in \mathbb{R}^{D};
\label{eq_KKT_unified}
\end{equation}%
\emph{where (\ref{eq_kkt2}) reduces to just }$\left\{ Qx=-c,Ax\leq b\right\} 
$\emph{\ in the case when }$D=\emptyset .$

$\left( ii\right) $\emph{\ Observe that any }$D\subset \{1,...,m\}$\emph{\
involved in (\ref{eq_kkt2}) satisfies} $D\subset I_{b}\left( x\right) $ 
\emph{since one of the blocks of (\ref{eq_kkt2}) is }$A_{D}x=b_{D}.$
\end{rem}

For each $D\subset I=\{1,...,m\},$ we consider the canonically perturbed
problem 
\begin{equation}
\begin{array}{cl}
\left( P_{D}\left( c,b_{D}\right) \right) & \text{minimize }\frac{1}{2}%
x^{\prime }Qx+c^{\prime }x \\ 
& \text{s.t. ~~~\ ~~~~~~}A_{D}x\leq b_{D},%
\end{array}
\label{eq_PD}
\end{equation}%
and the associated optimal set mapping $\mathcal{S}_{D}:\mathbb{R}^{n}\times 
\mathbb{R}^{D}\rightrightarrows \mathbb{R}^{n},$\emph{\ }assigning to each $%
\left( c,b_{D}\right) $ the set of all optimal solutions of problem $%
P_{D}\left( c,b_{D}\right) $. Note that, $P_{\emptyset }\left( c\right) $ is
nothing else but the unconstrained optimization problem given by minimizing
the function $x\mapsto \frac{1}{2}x^{\prime }Qx+c^{\prime }x$ on $\mathbb{R}%
^{n}.$ When $D=\{1,...,m\}$ we just write $P\left( c,b\right) .$ Observe
that an element of $\mathcal{S}_{D}\left( c,b_{D}\right) $ is not
necessarily a feasible solution of our original problem $P\left( c,b\right) $%
.

From now we fix $\left( \left( \overline{c},\overline{b}\right) ,\overline{x}%
\right) \in \mathrm{gph}\mathcal{S}.$ Our primary focus is on those $%
\mathcal{S}_{D}$ with $D\in \mathcal{M}_{\overline{c},\overline{b}}\left( 
\overline{x}\right) ;$ in particular, in such a case $A_{D}$ has full row
rank as a consequence of Lemma \ref{Lem1}$\left( i\right) .$ For didactical
reasons, when $D\in \mathcal{M}_{\overline{c},\overline{b}}\left( \overline{x%
}\right) ,$ problem $P_{D}\left( c,b_{D}\right) $ will be referred sometimes
to as the \emph{minimal KKT subproblem} associated with $D$ at parameter $%
\left( c,b_{D}\right) $ (as far as it is constrained by a subsystem
associated with a minimal KKT set of indices). Later on, in Section 5, we
will need to extend the family $\mathcal{M}_{\overline{c},\overline{b}%
}\left( \overline{x}\right) .$

\section{Aubin property for minimal KKT subproblems}

We start by characterizing the Aubin property of $\mathcal{S}_{D}$ at a
point $\left( \left( \overline{c},\overline{b}_{D}\right) ,\overline{x}%
\right) $ of its graph with $\left( \left( \overline{c},\overline{b}\right) ,%
\overline{x}\right) \in \mathrm{gph}\mathcal{S}$ and $D\in \mathcal{M}_{%
\overline{c},\overline{b}}\left( \overline{x}\right) .$ From now on, for
simplicity in the notation, we write 
\begin{equation*}
M_{D}:=\left( 
\begin{array}{cc}
Q & A_{D}^{\prime } \\ 
A_{D} & 0_{D\times D}%
\end{array}%
\right)
\end{equation*}%
for all $D\subset \{1,...,m\},$ with $M_{\emptyset }:=Q.$

\begin{theo}
\label{Theorem_SD}Let $\left( \left( \overline{c},\overline{b}\right) ,%
\overline{x}\right) \in \mathrm{gph}\mathcal{S}$ and consider any $D\in 
\mathcal{M}_{\overline{c},\overline{b}}\left( \overline{x}\right) .$ The
following conditions are equivalent:

$\left( i\right) $ $\mathcal{S}_{D}$ has the Aubin property at $\left(
\left( \overline{c},\overline{b}_{D}\right) ,\overline{x}\right) ;$

$\left( ii\right) $ $\mathcal{S}_{D}$ is strongly Lipschitz stable at $%
\left( \left( \overline{c},\overline{b}_{D}\right) ,\overline{x}\right) ;$

$\left( iii\right) $ $\mathcal{S}_{D}$ is single valued in some neighborhood
of $\left( \overline{c},\overline{b}_{D}\right) ;$

$\left( iv\right) $ $\mathcal{S}_{D}\left( \overline{c},\overline{b}%
_{D}\right) =\{\overline{x}\};$

$\left( v\right) $ $M_{D}$ is nonsingular.
\end{theo}

\begin{dem}
$\left( i\right) \Rightarrow \left( ii\right) $ was already established in
Proposition \ref{Prop1}, by just replacing $\mathcal{S}$ with $\mathcal{S}%
_{D}.\,\ $

$\left( ii\right) \Rightarrow \left( iii\right) $ and $\left( iii\right)
\Rightarrow \left( iv\right) $ are trivial.

Let us see $\left( iv\right) \Rightarrow \left( v\right) .$ In the case when 
$D=\emptyset $, 
\begin{equation*}
\mathcal{S}_{\emptyset }\left( \overline{c}\right) =\left\{ x\in \mathbb{R}%
^{n}\mid Qx=-\overline{c}\right\} ,
\end{equation*}%
and its single-valuedness entails that $M_{\emptyset }:=Q$ is nonsingular.

Assume now that $D\neq \emptyset .$ Let $\overline{\lambda }=\left( 
\overline{\lambda }_{i}\right) _{i\in D}$ be such that $\left( \overline{x},%
\overline{\lambda }\right) $ is a solution of (\ref{eq_kkt2}) with $\left(
c,b\right) =\left( \overline{c},\overline{b}\right) $ therein. The
minimality of $D,$ ensures $\overline{\lambda }_{i}>0$ for all $i\in D.$
Arguing by contradiction, assume that $M_{D}$ is singular, yielding the
existence of another solution $\left( \widetilde{x},\widetilde{\lambda }%
\right) \neq \left( \overline{x},\overline{\lambda }\right) $ of the linear
system in the variable $\left( x,\lambda \right) $ 
\begin{equation}
M_{D}\dbinom{x}{\lambda }=\dbinom{-\overline{c}}{\overline{b}_{D}}.
\label{eq_1}
\end{equation}%
For each $r\in \mathbb{N},$ the convex combination $\left( x^{r},\lambda
^{r}\right) :=\left( 1-\frac{1}{r}\right) \left( \overline{x},\overline{%
\lambda }\right) +\frac{1}{r}\left( \widetilde{x},\widetilde{\lambda }%
\right) \in \mathbb{R}^{n}\times \mathbb{R}^{D}$ gives a solution of (\ref%
{eq_1}).\thinspace\ For $r$ large enough, say $r\geq r_{0},$ we have that
all coordinates of $\lambda ^{r},$ say $\lambda _{i}^{r},$ $i\in D,$ are
positive. Observe that, for $r\geq r_{0},$ $x^{r}\neq \overline{x}$ since
otherwise we would have that $\left( \overline{x},\lambda ^{r}\right) $ is a
solution of (\ref{eq_1}) different from $\left( \overline{x},\overline{%
\lambda }\right) $, which contradicts that $A_{D}$ has full row rank (recall
that $D\in \mathcal{M}_{\overline{c},\overline{b}}\left( \overline{x}\right) 
$ and Lemma \ref{Lem1}).\ More in detail, we would have $A_{D}^{\prime
}\left( \lambda ^{r}-\overline{\lambda }\right) =0_{n}.$

Therefore, appealing to the KKT conditions for problem $P_{D}\left( 
\overline{c},\overline{b}_{D}\right) ,$ we have that 
\begin{equation*}
\overline{x},x^{r}\in \mathcal{S}_{D}\left( \overline{c},\overline{b}%
_{D}\right) ,\text{ for }r\text{ large enough,}
\end{equation*}%
which contradicts $\left( iv\right) .$

$\left( v\right) \Rightarrow \left( i\right) .$ Similarly to the previous
implication, the only solution to (\ref{eq_1}), say $\left( \overline{x},%
\overline{\lambda }\right) ,$ satisfies $\overline{\lambda }_{i}>0$ for all $%
i\in D.$ Let us consider $\left( c,b_{D}\right) $ close enough to $\left( 
\overline{c},\overline{b}_{D}\right) ,$ say $\left( c,b_{D}\right) $
belonging to a certain neighborhood $V$ of $\left( \overline{c},\overline{b}%
_{D}\right) ,$ to ensure that $\lambda _{D}\left( c,b_{D}\right) \geq 0_{D}$
whenever $\left( c,b_{D}\right) \in V,$ where 
\begin{equation}
\dbinom{x_{D}\left( c,b_{D}\right) }{\lambda _{D}\left( c,b_{D}\right) }%
:=M_{D}^{-1}\dbinom{-c}{b_{D}}.  \label{eq_3}
\end{equation}%
Clearly, $x_{D}\left( c,b_{D}\right) \in \mathcal{S}_{D}\left(
c,b_{D}\right) $ for $\left( c,b_{D}\right) \in V.$ The Aubin property of $%
\mathcal{S}_{D}$ at $\left( \left( \overline{c},\overline{b}_{D}\right) ,%
\overline{x}\right) $ follows immediately from (\ref{eq_3}) if we see that 
\begin{equation}
\mathcal{S}_{D}\left( c,b_{D}\right) =\left\{ x_{D}\left( c,b_{D}\right)
\right\}  \label{eq_4}
\end{equation}%
for $\left( c,b_{D}\right) $ in a possibly smaller neighborhood of $\left( 
\overline{c},\overline{b}_{D}\right) .$ To this end, we assume to the
contrary that there exist sequences $\{(c^{r},b_{D}^{r})\}_{r\in \mathbb{N}}$
and $\{x^{r}\}_{r\in \mathbb{N}}$ with $x^{r}\in \mathcal{S}_{D}\left(
c^{r},b_{D}^{r}\right) \backslash \{x_{D}\left( c^{r},b_{D}^{r}\right) \},$
and $\{(c^{r},b_{D}^{r})\}_{r\in \mathbb{N}}$ converging to $\left( 
\overline{c},\overline{b}_{D}\right) .$ For each $r\in \mathbb{N},$ the
convex combination 
\begin{equation*}
y^{r}:=\left( 1-\frac{1}{r\left( \left\Vert x^{r}\right\Vert +1\right) }%
\right) x_{D}\left( c^{r},b_{D}^{r}\right) +\frac{1}{r\left( \left\Vert
x^{r}\right\Vert +1\right) }x^{r}
\end{equation*}%
belongs to $\mathcal{S}_{D}\left( c^{r},b_{D}^{r}\right) ,$ since this set
is convex. Then, the KKT$\ $conditions applied to $P_{D}\left(
c^{r},b_{D}^{r}\right) $ ensure the existence of $D_{r}\subset D$ and $%
\lambda ^{r}=\left( \lambda _{i}^{r}\right) _{i\in D_{r},}\in \mathbb{R}%
_{+}^{D_{r}}$ such that 
\begin{equation}
M_{D_{r}}\dbinom{y^{r}}{\lambda ^{r}}=\dbinom{-c^{r}}{b_{D_{r}}^{r}},\text{
and }A_{D\backslash D_{r}}y^{r}\leq b_{D\backslash D_{r}}^{r}.  \label{eq_2}
\end{equation}%
In fact, 
\begin{equation*}
D_{r}\varsubsetneqq D,\text{ for all }r,
\end{equation*}%
since $y^{r}\neq x_{D}\left( c^{r},b_{D}^{r}\right) $ (the contrary would
contradict the definition of\linebreak $x_{D}\left( c^{r},b_{D}^{r}\right) $%
). The finiteness of $D$ allows us to assume that $\left\{ D_{r}\right\}
_{r\in \mathbb{N}}$ is constant (by taking a subsequence if necessary), say $%
D_{r}=\widetilde{D}\varsubsetneqq D$ for all $r.$ Let us see that the
sequence $\left\{ \lambda ^{r}\right\} _{r\in \mathbb{N}}$ is bounded.
Otherwise, we may assume without loss of generality $\sum_{i\in \widetilde{D}%
}\lambda _{i}^{r}\rightarrow \infty $ and $\left( \sum_{i\in \widetilde{D}%
}\lambda _{i}^{r}\right) ^{-1}\lambda _{j}^{r}\rightarrow \mu _{j}\geq 0$
for each $j\in \widetilde{D},$ with $\sum_{j\in \widetilde{D}}\mu _{j}=1.$
Dividing, for each $r,$ both sides of (\ref{eq_2})\thinspace by $\sum_{i\in 
\widetilde{D}}\lambda _{i}^{r}$ and letting $r\rightarrow \infty $ we obtain
the contradiction (recall that $A_{\widetilde{D}}$ has also full row rank,
as $A_{D}$) 
\begin{equation*}
A_{\widetilde{D}}^{\prime }\mu =0_{n},\text{ with }\mu =\left( \mu
_{j}\right) _{j\in \widetilde{D}}.
\end{equation*}%
Here we have used the boundedness of sequences $\left\{ y^{r}\right\} ,$ $%
\left\{ c^{r}\right\} ,$ and $\left\{ b_{D}^{r}\right\} .$ Observe that (\ref%
{eq_3}) entails $x_{D}\left( c^{r},b_{D}^{r}\right) \rightarrow x_{D}\left( 
\overline{c},\overline{b}_{D}\right) =\overline{x}$ as $r\rightarrow \infty
, $ and the definition of $y^{r}$ yields $\left\Vert y^{r}-x_{D}\left(
c^{r},b_{D}^{r}\right) \right\Vert \rightarrow 0.$ Once we know that $%
\left\{ \lambda ^{r}\right\} _{r\in \mathbb{N}}$ is bounded, we may assume $%
\lambda ^{r}\rightarrow \widetilde{\lambda }\in \mathbb{R}_{+}^{\widetilde{D}%
}$ and letting $r\rightarrow \infty $ in (\ref{eq_2}) we obtain%
\begin{equation*}
M_{\widetilde{D}}\dbinom{\overline{x}}{\widetilde{\lambda }}=\dbinom{-%
\overline{c}}{\overline{b}_{\widetilde{D}}}.
\end{equation*}%
This, together with the already known $A\overline{x}\leq \overline{b}$
yields a contradiction with the minimality assumption $D\in \mathcal{M}_{%
\overline{c},\overline{b}}\left( \overline{x}\right) .$
\end{dem}

Given $\left( \left( \overline{c},\overline{b}\right) ,\overline{x}\right)
\in \mathrm{gph}\mathcal{S}$ and $D\in \mathcal{M}_{\overline{c},\overline{b}%
}\left( \overline{x}\right) $ such that $\mathcal{S}_{D}$ has the Aubin
property at $\left( \left( \overline{c},\overline{b}_{D}\right) ,\overline{x}%
\right) ,$ the previous proof establishes the existence of a function 
\begin{equation}
\left( c,b\right) \mapsto \left( x_{D}\left( c,b_{D}\right) ,\lambda
_{D}\left( c,b_{D}\right) \right)  \label{eq_norms_in_products}
\end{equation}%
from $\left( \mathbb{R}^{n}\times \mathbb{R}^{m},\max \left\{ \left\Vert
\cdot \right\Vert _{\ast },\left\Vert \cdot \right\Vert _{\infty }\right\}
\right) $ to $(\mathbb{R}^{n}\times \mathbb{R}^{D},\left\Vert \cdot
\right\Vert +\left\Vert \cdot \right\Vert _{1})$ defined by (\ref{eq_3}) in
a certain neighborhood $V$ of $\left( \overline{c},\overline{b}\right) .$
Then, by projecting on the first $n$ coordinates, 
\begin{equation*}
\mathcal{S}_{D}\left( c,b_{D}\right) =\left\{ \left( 
\begin{array}{cc}
I_{n} & 0_{n\times \left\vert D\right\vert }%
\end{array}%
\right) M_{D}^{-1}\dbinom{-c}{b_{D}}\right\} .
\end{equation*}%
Now, the following result is immediate from the definition of matrix norm.

\begin{cor}
\label{Cor_lipSD}Let $\left( \left( \overline{c},\overline{b}\right) ,%
\overline{x}\right) \in \mathrm{gph}\mathcal{S}$ and $D\in \mathcal{M}_{%
\overline{c},\overline{b}}\left( \overline{x}\right) $ such that $\mathcal{S}%
_{D}$ has the Aubin property at $\left( \left( \overline{c},\overline{b}%
_{D}\right) ,\overline{x}\right) .$ Then, one has%
\begin{equation*}
\mathrm{lip}\mathcal{S}_{D}\left( \left( \overline{c},\overline{b}%
_{D}\right) ,\overline{x}\right) =\left\Vert \left( 
\begin{array}{cc}
I_{n} & 0_{n\times \left\vert D\right\vert }%
\end{array}%
\right) M_{D}^{-1}\right\Vert .
\end{equation*}
\end{cor}

Recall that, with our current choice of norms (\ref{eq_norms}) the previous
matrix norm is 
\begin{equation*}
\max_{\max \left\{ \left\Vert \alpha \right\Vert _{\ast },\left\Vert \beta
\right\Vert _{\infty }\right\} =1}\left\Vert \left( 
\begin{array}{cc}
I_{n} & 0_{n\times \left\vert D\right\vert }%
\end{array}%
\right) M_{D}^{-1}\dbinom{\alpha }{\beta }\right\Vert ,
\end{equation*}%
where $\left\Vert \cdot \right\Vert $ is the arbitrary norm being considered
in $\mathbb{R}^{n}$ as the space of decision variables $x.$

\begin{rem}
\emph{In the linear case, when }$Q=0_{n\times n},$ \emph{and under the NC,
the previous expression reduces to} 
\begin{equation*}
\mathrm{lip}\mathcal{S}_{D}\left( \left( \overline{c},\overline{b}%
_{D}\right) ,\overline{x}\right) =\left\Vert A_{D}^{-1}\right\Vert
\end{equation*}%
\emph{(see Theorem \ref{Th_CGP08} and recall that the NC entails that }$%
A_{D} $\emph{\ is invertible for all }$D\in \mathcal{M}_{\overline{c},%
\overline{b}}\left( \overline{x}\right) $\emph{).}

\emph{In fact, in such a case,} 
\begin{eqnarray*}
\mathrm{lip}\mathcal{S}_{D}\left( \left( \overline{c},\overline{b}%
_{D}\right) ,\overline{x}\right) &=&\left\Vert \left( 
\begin{array}{cc}
I_{n} & 0_{n\times \left\vert D\right\vert }%
\end{array}%
\right) \left( 
\begin{array}{cc}
0_{n\times n} & A_{D}^{-1} \\ 
\left( A_{D}^{-1}\right) ^{\prime } & 0_{D\times D}%
\end{array}%
\right) \right\Vert \\
&=&\max_{\max \left\{ \left\Vert \alpha \right\Vert _{\ast },\left\Vert
\beta \right\Vert _{\infty }\right\} =1}\left\Vert \left( 
\begin{array}{cc}
0_{n\times n} & A_{D}^{-1}%
\end{array}%
\right) \dbinom{\alpha }{\beta }\right\Vert \\
&=&\max_{\left\Vert \beta \right\Vert _{\infty }=1}\left\Vert
A_{D}^{-1}\beta \right\Vert =\left\Vert A_{D}^{-1}\right\Vert .
\end{eqnarray*}
\end{rem}

\begin{rem}
\label{Rem_S_SD}$\left( i\right) $ \emph{For an arbitrary} $D\subset
\{1,....,m\}$ \emph{we can have} $\mathcal{S}_{D}\left( c,b_{D}\right) \cap 
\mathcal{F}\left( b\right) =\emptyset ;$ \emph{for instance, consider the
problem, in} $\mathbb{R}^{2},$ \emph{of minimizing} $x_{1}^{2}+x_{2}^{2}$ 
\emph{subject to }$-x_{1}-x_{2}\leq -1,$ $-x_{1}\leq 0,$ $-x_{2}\leq 0$ 
\emph{and take }$D=\{2,3\}.$

$\left( ii\right) $ \emph{In the case when }$D\in \mathcal{M}_{c,b}\left(
x\right) ,$ \emph{for }$\left( \left( c,b\right) ,x\right) \in \limfunc{gph}%
\mathcal{S},$\emph{\ we always have }$x\in \mathcal{S}_{D}\left(
c,b_{D}\right) \cap \mathcal{F}\left( b\right) $ \emph{since }$A_{D}x=b_{D}$ 
\emph{and} $-(Qx+c)\in \mathrm{cone}\left\{ a_{i},i\in D\right\} $ \emph{%
yield to the KKT conditions for }$P_{D}\left( c,b_{D}\right) $ \emph{at }$x.$

$\left( iii\right) $ \emph{For }$D\subset \{1,....,m\},$ \emph{condition }$%
\mathcal{S}_{D}\left( c,b_{D}\right) \cap \mathcal{F}\left( b\right) \neq
\emptyset $ \emph{implies} $\mathcal{S}_{D}\left( c,b_{D}\right) \cap 
\mathcal{F}\left( b\right) =\mathcal{S}\left( c,b\right) .\,$\emph{In\ fact,
the inclusion `}$\subset $\emph{' is clear, while the converse comes from
the fact that both optimization problems have the same objective function
(and, hence, the same optimal value when }$\mathcal{S}_{D}\left(
c,b_{D}\right) \cap \mathcal{F}\left( b\right) \neq \emptyset $\emph{).}
\end{rem}

The following proposition includes the possibility that $\mathrm{lip}%
\mathcal{S}_{D}\left( \left( \overline{c},\overline{b}_{D}\right) ,\overline{%
x}\right) =+\infty $ for some $D\in \mathcal{M}_{\bar{c},\bar{b}}\left( \bar{%
x}\right) ,$ in which case the statement of the proposition yields $\mathrm{%
lip}\mathcal{S}\left( \left( \overline{c},\overline{b}\right) ,\overline{x}%
\right) =+\infty .$

\begin{prop}
\label{Prop>=}Let $\left( \left( \overline{c},\overline{b}\right) ,\overline{%
x}\right) \in \mathrm{gph}\mathcal{S}.$ One has 
\begin{equation}
\mathrm{lip}\mathcal{S}\left( \left( \overline{c},\overline{b}\right) ,%
\overline{x}\right) \geq \max_{D\in \mathcal{M}_{\bar{c},\bar{b}}\left( \bar{%
x}\right) }\mathrm{lip}\mathcal{S}_{D}\left( \left( \overline{c},\overline{b}%
_{D}\right) ,\overline{x}\right) .  \label{eq_5}
\end{equation}
\end{prop}

\begin{dem}
First, assume that $\mathrm{lip}\mathcal{S}_{D}\left( \left( \overline{c},%
\overline{b}_{D}\right) ,\overline{x}\right) =+\infty \,\ $for a certain $%
D\in \mathcal{M}_{\bar{c},\bar{b}}\left( \bar{x}\right) ,$ i.e., that $%
\mathcal{S}_{D}$ has not the Aubin property at $\left( \left( \overline{c},%
\overline{b}_{D}\right) ,\overline{x}\right) ,$ and let us see that $%
\mathcal{S}$ also has not the Aubin property at $\left( \left( \overline{c},%
\overline{b}\right) ,\overline{x}\right) .$ Under the current assumption,
Theorem \ref{Theorem_SD} $\left( iv\right) $ ensures the existence of $%
\widetilde{x}\in \mathcal{S}_{D}\left( \overline{c},\overline{b}_{D}\right) $
with $\widetilde{x}\neq \overline{x}.$ Consider, for each $r\in \mathbb{N},$ 
$x^{r}:=\left( 1-\frac{1}{r}\right) \overline{x}+\frac{1}{r}\widetilde{x}\in 
\mathcal{S}_{D}\left( \overline{c},\overline{b}_{D}\right) $ and $%
b^{r}=\left( b_{j}^{r}\right) _{j=1,...,m}\in \mathbb{R}^{m}$ defined as%
\begin{equation*}
b_{j}^{r}:=\left\{ 
\begin{array}{l}
\overline{b}_{j},\text{ if }j\in D, \\ 
\overline{b}_{j}+\left[ a_{j}^{\prime }x^{r}-\overline{b}_{j}\right] _{+},%
\text{ if }j\in \{1,...,m\}\backslash D,%
\end{array}%
\right.
\end{equation*}%
where $\left[ \gamma \right] _{+}:=\max \{\gamma ,0\}.$ In this way, for
each $r,$ $x^{r}\in \mathcal{F}\left( b^{r}\right) ;$ moreover $x^{r}\in 
\mathcal{S}_{D}\left( \overline{c},b_{D}^{r}\right) $ since $b_{D}^{r}=%
\overline{b}_{D}.\,\ $So, $x^{r}\in \mathcal{S}_{D}\left( \overline{c}%
,b_{D}^{r}\right) \cap \mathcal{F}\left( b^{r}\right) \subset \mathcal{S}%
\left( \overline{c},b^{r}\right) ;$ see $\left( iii\right) $ in the previous
remark. On the other hand, we also have $\overline{x}\in \mathcal{S}%
_{D}\left( \overline{c},b_{D}^{r}\right) \cap \mathcal{F}\left( b^{r}\right)
,$ and then $\overline{x}\in \mathcal{S}\left( \overline{c},b^{r}\right) ,$
for each $r.$ Since $\left\{ b^{r}\right\} $ converges to $\overline{b}$ (as
far as $\left\{ x^{r}\right\} $ converges to $\overline{x}),$ $\mathcal{S}$
is not single-valued in any neighborhood of $\left( \overline{c},\overline{b}%
\right) ,$ and then $\mathcal{S}$ has not the Aubin property at $\left(
\left( \overline{c},\overline{b}\right) ,\overline{x}\right) $ (recall
Proposition \ref{Prop1}).

Now, we deal with the case $\mathrm{lip}\mathcal{S}_{D}\left( \left( 
\overline{c},\overline{b}_{D}\right) ,\overline{x}\right) <+\infty $ for all 
$D\in \mathcal{M}_{\bar{c},\bar{b}}\left( \bar{x}\right) .$ Fix any $D\in 
\mathcal{M}_{\overline{c},\overline{b}}\left( \overline{x}\right) $ and
appealing to Theorem \ref{Theorem_SD}$\left( iii\right) ,$ write 
\begin{equation*}
\mathrm{lip}\mathcal{S}_{D}\left( \left( \overline{c},\overline{b}%
_{D}\right) ,\overline{x}\right) =\lim_{r\rightarrow \infty }\frac{%
\left\Vert x^{r}-\widetilde{x}^{r}\right\Vert }{\left\Vert \left(
c^{r},\beta ^{r}\right) -\left( \widetilde{c}^{r},\widetilde{\beta }%
^{r}\right) \right\Vert },
\end{equation*}%
for some sequences of parameters $\left\{ \left( c^{r},\beta ^{r}\right)
\right\} ,$ $\left\{ \left( \widetilde{c}^{r},\widetilde{\beta }^{r}\right)
\right\} \subset \mathbb{R}^{n}\times \mathbb{R}^{D}$ converging to $\left( 
\overline{c},\overline{b}_{D}\right) ,$ such that $\mathcal{S}_{D}\left(
c^{r},\beta ^{r}\right) =\{x^{r}\},$ $\mathcal{S}_{D}\left( \widetilde{c}%
^{r},\widetilde{\beta }^{r}\right) =\{\widetilde{x}^{r}\},$ for all $r.$

For each $r,$ define $b^{r}=\left( b_{j}^{r}\right) _{j=1,...,m},$ $%
\widetilde{b}^{r}=\left( \widetilde{b}_{j}^{r}\right) _{j=1,...,m}\in 
\mathbb{R}^{m}$ as follows:%
\begin{equation*}
b_{j}^{r}:=\beta _{j}^{r}\text{ and }\widetilde{b}_{j}^{r}:=\widetilde{\beta 
}_{j}^{r},\text{ whenever }j\in D
\end{equation*}%
and%
\begin{equation*}
b_{j}^{r}=\widetilde{b}_{j}^{r}:=\overline{b}_{j}+\max \left\{ \left[
a_{j}^{\prime }x^{r}-\overline{b}_{j}\right] _{+},\left[ a_{j}^{\prime }%
\widetilde{x}^{r}-\overline{b}_{j}\right] _{+}\right\} ,\text{ for }j\in
\{1,...,m\}\backslash D.
\end{equation*}

In this way, one easily checks that $x^{r}\in \mathcal{F}\left( b^{r}\right) 
$ and $\widetilde{x}^{r}\in \mathcal{F}\left( \widetilde{b}^{r}\right) $ for
all $r.$ Moreover, for each $r,$ $x^{r}\in \mathcal{S}_{D}\left(
c^{r},b_{D}^{r}\right) $ and $\widetilde{x}^{r}\in \mathcal{S}_{D}\left(
c^{r},\widetilde{b}_{D}^{r}\right) ,$ since $b_{D}^{r}=\beta ^{r}$ and $%
\widetilde{b}_{D}^{r}=\widetilde{\beta }^{r}.$ Therefore, appealing to
Remark \ref{Rem_S_SD}$\left( iii\right) $ we have 
\begin{equation*}
\mathcal{S}\left( c^{r},b^{r}\right) =\{x^{r}\}\text{ and }\mathcal{S}\left( 
\widetilde{c}^{r},\widetilde{b}^{r}\right) =\{\widetilde{x}^{r}\}.
\end{equation*}%
Finally, observe that, for each $r,$%
\begin{equation*}
\left\Vert \left( c^{r},b^{r}\right) -\left( \widetilde{c}^{r},\widetilde{b}%
^{r}\right) \right\Vert =\left\Vert \left( c^{r},\beta ^{r}\right) -\left( 
\widetilde{c}^{r},\widetilde{\beta }^{r}\right) \right\Vert ,
\end{equation*}%
since $\left\Vert b^{r}-\widetilde{b}^{r}\right\Vert _{\infty }=\left\Vert
b_{D}^{r}-\widetilde{b}_{D}^{r}\right\Vert _{\infty }=\left\Vert \beta ^{r}-%
\widetilde{\beta }^{r}\right\Vert _{\infty },$ and consequently,%
\begin{equation*}
\mathrm{lip}\mathcal{S}\left( \left( \overline{c},\overline{b}\right) ,%
\overline{x}\right) \geq \lim_{r}\frac{\left\Vert x^{r}-\widetilde{x}%
^{r}\right\Vert }{\left\Vert \left( c^{r},b^{r}\right) -\left( \widetilde{c}%
^{r},\widetilde{b}^{r}\right) \right\Vert }=\mathrm{lip}\mathcal{S}%
_{D}\left( \left( \overline{c},\overline{b}_{D}\right) ,\overline{x}\right) .
\end{equation*}
\end{dem}

\begin{cor}
\label{CorS pL implies D nonsingular}Let $\left( \left( \overline{c},%
\overline{b}\right) ,\overline{x}\right) \in \mathrm{gph}\mathcal{S}.$ If $%
\mathcal{S}$ has the Aubin property at $\left( \left( \overline{c},\overline{%
b}\right) ,\overline{x}\right) $ then $M_{D}$ is nonsingular for all $D\in 
\mathcal{M}_{\overline{c},\overline{b}}\left( \overline{x}\right) .$
\end{cor}

\begin{dem}
As an immediate consequence of the previous proposition, if $\mathcal{S}$
has the Aubin property at $\left( \left( \overline{c},\overline{b}\right) ,%
\overline{x}\right) ,$ then for any $D\in \mathcal{M}_{\overline{c},%
\overline{b}}\left( \overline{x}\right) ,$ $\mathcal{S}_{D}$ also has the
Aubin property at $\left( \left( \overline{c},\overline{b}_{D}\right) ,%
\overline{x}\right) .$ Hence, \ref{Theorem_SD}$\left( v\right) $ yields the
aimed conclusion.
\end{dem}

The inequality in Proposition \ref{Prop>=} can be strict as the following
example shows.

\begin{exa}
\label{Exa1}\emph{Consider the following quadratic problem, in }$\mathbb{R}%
^{2}$ \emph{endowed with the Euclidean norm,}%
\begin{equation*}
\begin{tabular}{lll}
$P\left( c,b\right) :$ & \textrm{Minimize} & $\frac{1}{2}x^{\prime
}Qx+c^{\prime }x$ \\ 
& \textrm{subject to} & $Ax\leq b,$%
\end{tabular}%
\end{equation*}%
\emph{where }$Q=\left( 
\begin{array}{cc}
1 & 0 \\ 
0 & 1%
\end{array}%
\right) ,$ $A=\left( 
\begin{array}{cc}
-1 & 0 \\ 
0 & \frac{-1}{10}%
\end{array}%
\right) ,$ \emph{with nominal parameters\ }$\overline{c}=0_{2}$ \emph{and} $%
\overline{b}=\binom{-1}{0}.$ \emph{One easily sees that}$\,\mathcal{S}\left( 
\overline{c},\overline{b}\right) =\left\{ \overline{x}\right\} $\emph{\ with}%
\begin{equation*}
\overline{x}=\binom{1}{0}\text{\emph{\ and }}\mathcal{M}_{\bar{c},\bar{b}%
}\left( \bar{x}\right) =\left\{ 1\right\} .
\end{equation*}%
\emph{So, the right-hand side member of (\ref{eq_5}) equals}%
\begin{eqnarray*}
\mathrm{lip}\mathcal{S}_{\left\{ 1\right\} }\left( \left( \overline{c},%
\overline{b}_{\{1\}}\right) ,\overline{x}\right) &=&\left\Vert \left( 
\begin{array}{ccc}
1 & 0 & 0 \\ 
0 & 1 & 0%
\end{array}%
\right) \left( 
\begin{array}{ccc}
1 & 0 & -1 \\ 
0 & 1 & 0 \\ 
-1 & 0 & 0%
\end{array}%
\right) ^{-1}\right\Vert \\
&=&\left\Vert \left( 
\begin{array}{ccc}
0 & 0 & -1 \\ 
0 & 1 & 0%
\end{array}%
\right) \right\Vert =\sqrt{2}.
\end{eqnarray*}%
\emph{Recall that in the space of parameters we are considering the norm }$%
\left\Vert \left( c,b_{1}\right) \right\Vert =\max \left\{ \left\Vert
c\right\Vert _{\ast },\left\vert b_{1}\right\vert \right\} ,$ \emph{so that }%
\begin{equation*}
\left\Vert \left( 
\begin{array}{ccc}
0 & 0 & -1 \\ 
0 & 1 & 0%
\end{array}%
\right) \right\Vert =\sup \left\{ \sqrt{c_{2}^{2}+b_{1}^{2}}%
:c_{1}^{2}+c_{2}^{2}\leq 1,\text{ }\left\vert b_{1}\right\vert \leq
1\right\} .
\end{equation*}%
\emph{However, let us see that }$\mathrm{lip}\mathcal{S}\left( \left( 
\overline{c},\overline{b}\right) ,\overline{x}\right) >\sqrt{2}.$ \emph{%
Consider the following sequence of parameters }$\left\{ \left(
c^{r},b^{r}\right) \right\} _{r\in \mathbb{N}}$ \emph{given by} 
\begin{equation*}
\left( c^{r},b^{r}\right) =\left( 0_{2},\left( -1-\frac{1}{r},-\frac{1}{r}%
,\right) ^{\prime }\right) ,\text{ \emph{for each} }r.
\end{equation*}%
\emph{One easily sees that }$\mathcal{S}\left( c^{r},b^{r}\right) =\left\{
\left( 1+\frac{1}{r},\frac{10}{r}\right) ^{\prime }\right\} ,$ \emph{for all}
$r,$ \emph{and then}%
\begin{equation*}
\mathrm{lip}\mathcal{S}\left( \left( \overline{c},\overline{b}\right) ,%
\overline{x}\right) \geq \lim \frac{\left\Vert \left( 1+\frac{1}{r},\frac{10%
}{r}\right) ^{\prime }-\left( 1,0\right) ^{\prime }\right\Vert }{\frac{1}{r}}%
=\sqrt{101}.
\end{equation*}%
\emph{In fact, one can check by applying Theorem \ref{Th main} (in Section
6) that the previous inequality is indeed an equality.}
\end{exa}

\section{Extended subproblems}

Example \ref{Exa1} shows that the family $\mathcal{M}_{\overline{c},%
\overline{b}}\left( \overline{x}\right) $ is not enough to
determine\linebreak $\mathrm{lip}\mathcal{S}\left( \left( \overline{c},%
\overline{b}\right) ,\overline{x}\right) $. A natural idea in the pursuit of
an exact formula consists of considering the maximum of $\mathrm{lip}%
\mathcal{S}_{D}\left( \left( \overline{c},\overline{b}_{D}\right) ,\overline{%
x}\right) $ over an extended family of subsets $D.$ Inspired by what happens
in this example, we consider the following extended family of KKT subsets of
indices 
\begin{equation*}
\mathcal{L}_{\overline{c},\overline{b}}\left( \overline{x}\right) :=\left\{
D\subset I_{\overline{b}}\left( \overline{x}\right) \left\vert 
\begin{array}{l}
\left\{ a_{i},i\in D\right\} \text{ is linearly independent} \\ 
\text{and }-(Q\overline{x}+\overline{c})\in \mathrm{cone}\left\{ a_{i},i\in
D\right\}%
\end{array}%
\right. \right\} .
\end{equation*}%
Directly from the definitions, we have%
\begin{equation}
\mathcal{M}_{\overline{c},\overline{b}}\left( \overline{x}\right) \subset 
\mathcal{L}_{\overline{c},\overline{b}}\left( \overline{x}\right) .
\label{eq_M in L}
\end{equation}%
Observe that, when dealing with any $D\in \mathcal{L}_{\overline{c},%
\overline{b}}\left( \overline{x}\right) ,$ the (unique) coefficients
generating $-(Q\overline{x}+\overline{c})$ as a linear combination of $%
\left\{ a_{i},i\in D\right\} $ are nonnegative, but some of them could be
zero.

In the following result we refer to the (sequential) Painlev\'{e}-Kuratowski
upper/outer limit of sets (see, e.g., \cite[p. 13]{mor06a}). For a generic
set-valued mapping $\mathcal{G}:Y\rightrightarrows X$ between metric spaces, 
$\limfunc{Limsup}_{y\rightarrow \overline{y}}\mathcal{G}\left( y\right) $
consists of those $x\in X$ for which there exist sequences $\left\{ \left(
y^{r},x^{r}\right) \right\} _{r\in \mathbb{N}}\subset \limfunc{gph}\mathcal{G%
}$ converging to $\left( \overline{y},x\right) .$ We do not exclude constant
sequences, hence $\mathcal{G}\left( \overline{y}\right) \subset \limfunc{%
Limsup}_{y\rightarrow \overline{y}}\mathcal{G}\left( y\right) .$

\begin{prop}
\label{Prop Limsup minimimals}Let $\left( \left( \overline{c},\overline{b}%
\right) ,\overline{x}\right) \in \mathrm{gph}\mathcal{S}.$ Then

$\left( i\right) $ There exists a neighborhood $V\times U$ of $\left( \left( 
\overline{c},\overline{b}\right) ,\overline{x}\right) $ such that 
\begin{equation}
\mathcal{M}_{c,b}\left( x\right) \subset \mathcal{L}_{\overline{c},\overline{%
b}}\left( \overline{x}\right) \text{ for all }\left( \left( c,b\right)
,x\right) \in \left( V\times U\right) \cap \mathrm{gph}\mathcal{S}.
\label{eq_MinL}
\end{equation}

$\left( ii\right) $%
\begin{equation*}
\mathcal{L}_{\overline{c},\overline{b}}\left( \overline{x}\right) =\limfunc{%
Limsup}_{\substack{ \left( \left( c,b\right) ,x\right) \rightarrow \left(
\left( \overline{c},\overline{b}\right) ,\overline{x}\right)  \\ \left(
\left( c,b\right) ,x\right) \in \mathrm{gph}\mathcal{S}}}\mathcal{M}%
_{c,b}\left( x\right) ,
\end{equation*}%
where the $\limfunc{Limsup}$\ is taken in the discrete space $2^{\left\{
1,...,m\right\} }.$

$\left( iii\right) $ For any open neighborhood $V\times U$ of $\left( \left( 
\overline{c},\overline{b}\right) ,\overline{x}\right) $ such that (\ref%
{eq_MinL}) holds, we have 
\begin{equation*}
\mathcal{L}_{c,b}\left( x\right) \subset \mathcal{L}_{\overline{c},\overline{%
b}}\left( \overline{x}\right) \text{ for all }\left( \left( c,b\right)
,x\right) \in \left( V\times U\right) \cap \mathrm{gph}\mathcal{S}.
\end{equation*}
\end{prop}

\begin{dem}
$\left( i\right) $ Reasoning by contradiction, assume that there exists a
sequence $\left\{ \left( \left( c^{r},b^{r}\right) ,x^{r}\right) \right\}
_{r\in \mathbb{N}}$ converging to $\left( \left( \overline{c},\overline{b}%
\right) ,\overline{x}\right) $ where, for all $r,$ we have $\left( \left(
c^{r},b^{r}\right) ,x^{r}\right) \in \mathrm{gph}\mathcal{S}$ and there
exists $D_{r}\in \mathcal{M}_{c^{r},b^{r}}\left( x^{r}\right) \backslash 
\mathcal{L}_{\overline{c},\overline{b}}\left( \overline{x}\right) .$ The
finiteness of $I=\{1,...,m\}$ allows us to assume (by taking a subsequence
if necessary) that $\left\{ D_{r}\right\} _{r\in \mathbb{N}}$ is constant,
say $D_{r}=D$ for all $r.$ From $D\in I_{b^{r}}\left( x^{r}\right) ,$ for
all $r,$ we easily obtain $D\subset I_{\overline{b}}\left( \overline{x}%
\right) .$ Let us show that $D\in \mathcal{L}_{\overline{c},\overline{b}%
}\left( \overline{x}\right) ,$ which constitutes a contradiction. To do
this, write 
\begin{equation}
-\left( Qx^{r}+c^{r}\right) =\sum_{i\in D}\lambda _{i}^{r}a_{i}  \label{eq10}
\end{equation}%
for some sequence $\left\{ \lambda ^{r}\right\} _{r\in \mathbb{N}}\subset 
\mathbb{R}_{+}^{D}$, where $\mathbb{R}_{+}:=\left[ 0,+\infty \right[ .$
Write $\gamma ^{r}:=\sum_{i\in D}\lambda _{i}^{r}$ for all $r\in \mathbb{N}.$
If $\left\{ \gamma ^{r}\right\} _{r\in \mathbb{N}}$ is unbounded in $\mathbb{%
R},$ we may assume $\gamma ^{r}\rightarrow +\infty $ by taking a suitable
subsequence, and dividing both members of (\ref{eq10}) by $\gamma ^{r}$ and
letting $r\rightarrow \infty $ we get $0_{n}\in \mathrm{conv}\left\{
a_{i},~i\in D\right\} ,$ contradicting our general assumption made after
Proposition \ref{Prop_SCQ}. Hence we may assume without loss of generality
that $\left\{ \lambda _{i}^{r}\right\} _{r\in \mathbb{N}}$ converges to some 
$\mu _{i}\geq 0$ for each $i\in D.$ Now (\ref{eq10}) gives 
\begin{equation*}
-\left( Q\overline{x}+\overline{c}\right) =\sum_{i\in D}\mu _{i}a_{i}.
\end{equation*}%
The linear independence of $\left\{ a_{i},i\in D\right\} $ comes from $D\in 
\mathcal{M}_{c^{r},b^{r}}\left( x^{r}\right) $ for any $r.$ Then we have $%
D\in \mathcal{L}_{\overline{c},\overline{b}}\left( \overline{x}\right) .$
Note that we cannot ensure that $D$ belongs to $\mathcal{M}_{\overline{c},%
\overline{b}}\left( \overline{x}\right) $ since some $\mu _{i}$ could be
zero.

$\left( ii\right) $ Inclusion `$\supset $' comes directly from $\left(
i\right) $ and the fact that the only convergent sequences in $2^{\left\{
1,...,m\right\} }$ are the constant ones. Reciprocally, given any $D\in 
\mathcal{L}_{\overline{c},\overline{b}}\left( \overline{x}\right) ,$ let us
write $-(Q\overline{x}+\overline{c})=\sum_{i\in D}\lambda _{i}a_{i}$ for
some nonnegative scalars $\left\{ \lambda _{i}\right\} _{i\in D}$. Then,
because of the linear independence of $\left\{ a_{i},i\in D\right\} ,$ it is
clear that $\left\{ i\in D:\lambda _{i}\neq 0\right\} \in \mathcal{M}_{%
\overline{c},\overline{b}}\left( \overline{x}\right) .$ By defining $c^{r}:=%
\overline{c}-\sum_{i\in D,\lambda _{i}=0}\frac{1}{r}a_{i},$ we can easily
check that $D\in \mathcal{M}_{c^{r},\overline{b}}\left( \overline{x}\right) $
for all $r\in \mathbb{N}.$ This proves `$\subset $'.

$\left( iii\right) $ Given any $\left( \left( c,b\right) ,x\right) \in
\left( V\times U\right) \cap \mathrm{gph}\mathcal{S}$ and any $D\in \mathcal{%
L}_{c,b}\left( x\right) ,$ following the same argument as in $\left(
ii\right) $ we can find $\widetilde{c}$ close enough to $c$ to ensure $%
\left( \left( \widetilde{c},b\right) ,x\right) \in \left( V\times U\right)
\cap \mathrm{gph}\mathcal{S}$ and $D\in \mathcal{M}_{\widetilde{c},b}\left(
x\right) .$ Accordingly, $\left( i\right) $ entails $D\in \mathcal{L}_{%
\overline{c},\overline{b}}\left( \overline{x}\right) .$
\end{dem}

The following results are devoted to extend some of the properties already
established for $\mathcal{M}_{\overline{c},\overline{b}}\left( \overline{x}%
\right) $ to the new family $\mathcal{L}_{\overline{c},\overline{b}}\left( 
\overline{x}\right) .$

\begin{prop}
\label{Prop M_D}Let $D\subset I$ with $\left\{ a_{i},~i\in D\right\} \,$%
being linearly independent and write $d=\left\vert D\right\vert $ (the
cardinality of $D$). For all $x\in \mathbb{R}^{n}$ and all $y\in \mathbb{R}%
^{d}$ one has 
\begin{equation*}
M_{D}\dbinom{x}{y}=0_{n+d}\Leftrightarrow Qx=0_{n},\text{ }A_{D}x=0_{d},%
\text{ and }y=0_{d}.
\end{equation*}%
Consequently, $M_{D}$ is nonsingular if and only if $\ker Q\cap \ker
A_{D}=\{0_{n}\}.$
\end{prop}

\begin{dem}
The `$\Leftarrow $' implication is trivial. Recalling that $Q$ is assumed to
be symmetric and positive semidefinite, it has a (not necessarily unique)
Cholesky decomposition $LL^{\prime }$. Hence, if $M_{D}\dbinom{x}{y}%
=0_{n+d}, $ then one has $A_{D}x=0_{d}$ and 
\begin{eqnarray*}
0 &=&x^{\prime }0_{n}=x^{\prime }\left( Qx+A_{D}^{\prime }y\right)
=x^{\prime }Qx+\left( A_{D}x\right) ^{\prime }y=x^{\prime }LL^{\prime
}x+0_{n}^{\prime }y \\
&=&\left( L^{\prime }x\right) ^{\prime }L^{\prime }x\Rightarrow L^{\prime
}x=0_{n}\Rightarrow Qx=0_{n}\Rightarrow A_{D}^{\prime }y=0_{n}.
\end{eqnarray*}%
Since $A_{D}^{\prime }y$ is a linear combination of the columns of $%
A_{D}^{\prime },$ which are linearly independent, we obtain $y=0_{d}.$
\end{dem}

\begin{cor}
\label{Cor_D0}Let $D_{0}\subset D\subset I$ with $\left\{ a_{i},~i\in
D\right\} \,$being linearly independent. If $M_{D_{0}}$ is nonsingular, then 
$M_{D}$ is also nonsingular.
\end{cor}

\begin{dem}
Just observe that, denoting $d=\left\vert D\right\vert $ and $%
d_{0}=\left\vert D_{0}\right\vert ,$ 
\begin{equation}
M_{D}=\left( 
\begin{array}{ccc}
Q & A_{D_{0}}^{\prime } & A_{D\backslash D_{0}}^{\prime } \\ 
A_{D_{0}} & 0_{d_{0}\times d_{0}} & 0_{d_{0}\times \left( d-d_{0}\right) }
\\ 
A_{D\backslash D_{0}} & 0_{\left( d-d_{0}\right) \times d_{0}} & 0_{\left(
d-d_{0}\right) \times \left( d-d_{0}\right) }%
\end{array}%
\right) .  \label{eq_MD}
\end{equation}%
Writing $y\in \mathbb{R}^{d}$ as $\dbinom{y_{0}}{y_{1}},$ with $y_{0}\in 
\mathbb{R}^{d_{0}}$ and $y_{1}\in \mathbb{R}^{d-d_{0}},$ we have from
Proposition \ref{Prop M_D} that 
\begin{equation*}
M_{D}\dbinom{x}{y}=0_{n+d}\Rightarrow Qx=0_{n},\text{ }A_{D}x=0_{d},\text{
and }y=0_{d}.
\end{equation*}%
In particular, $A_{D_{0}}x=0_{d_{0}}\,\ $and $y_{0}=0_{d_{0}}.$ Hence,
applying again Proposition \ref{Prop M_D}, one has $M_{D_{0}}\dbinom{x}{y_{0}%
}=0_{n+d_{0}}.$ Since $M_{D_{0}}$ is assumed to be nonsingular, $\dbinom{x}{%
y_{0}}=0_{n+d_{0}}.$ Now, looking at the first row in equation $M_{D}\dbinom{%
x}{y}=0_{n+d},$ and taking (\ref{eq_MD}) into account, we get 
\begin{equation*}
Qx+A_{D_{0}}^{\prime }y_{0}+A_{D\backslash D_{0}}^{\prime
}y_{1}=A_{D\backslash D_{0}}^{\prime }y_{1}=0_{n},
\end{equation*}%
which implies $y_{1}=0_{d-d_{0}}$ because of the linear independence of $%
\left\{ a_{i},~i\in D\backslash D_{0}\right\} .$ Accordingly, $\dbinom{x}{y}%
=0_{n+d}$, which proves that $\ker M_{D}=\left\{ 0_{n+d}\right\} $ and,
therefore, $M_{D}$ is nonsingular.
\end{dem}

\begin{rem}
$\left( i\right) $ \emph{If }$D_{0}=\emptyset $\emph{, the previous
corollary reads as: if }$Q$ \emph{is nonsingular (hence positive definite),
then }$M_{D}$ \emph{is also nonsingular (but not necessarily positive
semidefinite) whenever }$A_{D}$ \emph{has full-row-rank.}

$\left( ii\right) $ \emph{If }$Q=0_{n\times n},$ \emph{the last part of
Proposition \ref{Prop M_D} reads as (provided that }$A_{D}$ \emph{has
full-row-rank):} 
\begin{equation*}
\left( 
\begin{array}{cc}
0_{n\times n} & A_{D}^{\prime } \\ 
A_{D} & 0_{D\times D}%
\end{array}%
\right) \text{ \emph{is nonsingular}}\Leftrightarrow A_{D}\text{ \emph{is
square and nonsingular.}}
\end{equation*}%
\emph{This, together with Theorem \ref{Theorem_SD}, constitutes an
alternative proof of the last part of Proposition \ref{Prop NC linear}.}

$\left( iii\right) $\emph{\ In the previous corollary we cannot switch the
roles of} $D_{0}$\emph{\ and }$D.$ \emph{For instance, consider the problem,
in} $\mathbb{R}^{2},$ \emph{of minimizing }$\frac{1}{2}x_{1}^{2}+x_{2}$ 
\emph{subject to }$-x_{2}\leq 0$\emph{, whose unique optimal solution is }$%
\overline{x}=0_{2}.$\emph{\ Then }$M_{\emptyset }=Q$ \emph{is singular
whereas} 
\begin{equation*}
M_{\{1\}}=\left( 
\begin{array}{ccc}
1 & 0 & 0 \\ 
0 & 0 & -1 \\ 
0 & -1 & 0%
\end{array}%
\right)
\end{equation*}%
\emph{is nonsingular. Observe that in this case} $\emptyset \notin \mathcal{M%
}_{\bar{c},\bar{b}}\left( \bar{x}\right) .$
\end{rem}

The next example shows that $\mathcal{S}_{D}$ may not have the Aubin
property at $\left( \left( \overline{c},\overline{b}_{D}\right) ,\overline{x}%
\right) \ $with $M_{D}$ being nonsingular.

\begin{exa}
\label{Exa 3/2}\emph{Consider the problem in} $\mathbb{R}^{2}$\emph{\ of
minimizing }$\frac{1}{2}x_{1}^{2}+x_{1}$\emph{\ subject to }$-x_{1}\leq 0$%
\emph{\ and }$x_{2}\leq 0.$ \emph{Here }$\mathcal{S}\left( \overline{c},%
\overline{b}\right) =\{0\}\times \left] -\infty ,0\right] $\emph{\ and
Proposition \ref{Prop1} ensures that }$\mathcal{S}=\mathcal{S}_{\{1,2\}}$%
\emph{\ has not the Aubin property at }$\mathcal{S}\left( \left( \overline{c}%
,\overline{b}\right) ,\overline{x}\right) $\emph{\ for any }$\overline{x}\in 
\mathcal{S}\left( \overline{c},\overline{b}\right) .$ \emph{However,} 
\begin{equation*}
M_{\{1,2\}}=\left( 
\begin{array}{cccc}
1 & 0 & -1 & 0 \\ 
0 & 0 & 0 & 1 \\ 
-1 & 0 & 0 & 0 \\ 
0 & 1 & 0 & 0%
\end{array}%
\right)
\end{equation*}%
\emph{is nonsingular. Observe that, for any }$\overline{x}\in \mathcal{S}%
\left( \overline{c},\overline{b}\right) ,$\emph{\ }$\mathcal{M}_{\overline{c}%
,\overline{b}}\left( \overline{x}\right) =\left\{ \{1\}\right\} ,$ $\mathcal{%
L}_{\overline{c},\overline{b}}\left( \overline{x}\right) =\left\{
\{1\},\{1,2\}\right\} $\emph{\ and }$M_{\{1\}}$\emph{\ is singular.}
\end{exa}

\section{Lipschitz modulus of the argmin mapping}

In the following paragraphs we introduce some concepts and results needed
for establishing our main theorem at the end of this section.

It is clear that, for any $D\subset I=\left\{ 1,...,m\right\} $ (including $%
D=\emptyset $ as in Remark \ref{Remark_empty}), the (possibly empty) set of
solutions $\left( \left( c,b\right) ,x,\lambda \right) \in \left( \mathbb{R}%
^{n}\times \mathbb{R}^{m}\right) \times \mathbb{R}^{n}\times \mathbb{R}^{m}$
to system (\ref{eq_kkt2}) is a polyhedral convex set. According to \cite[%
Theorem 19.3]{Rock70}, the projection, $P_{D},$ of such a solution set on $%
\left( \mathbb{R}^{n}\times \mathbb{R}^{m}\right) \times \mathbb{R}^{n}$ is
also a polyhedral convex set. Since the KKT optimality conditions yield%
\begin{equation*}
\limfunc{gph}\mathcal{S}=\bigcup\nolimits_{D\subset I}P_{D},
\end{equation*}%
$\limfunc{gph}\mathcal{S}$ is a so-called \emph{polyhedral set}, i.e., a
finite union of polyhedral convex sets. A classical result by Robinson \cite[%
Proposition 1]{Robinson} establishes that any \emph{polyhedral multifunction}%
, i.e., whose graph is a polyhedral set in $\mathbb{R}^{k},$ $k\in \mathbb{N}%
,$ is \emph{Lipschitz upper semicontinuous} at any point of its domain with
a common constant (the terminology used in \cite{Robinson} is \emph{locally
upper Lipschitzian}). Although a general definition can be given for
multifunctions between metric spaces, here we give the definition for our
mapping $\mathcal{S}$. We say that $\mathcal{S}$ is Lipschitz upper
semicontinuous at $\left( c,b\right) \in \limfunc{dom}\mathcal{S}$ with
constant $\kappa $ if there exists a neighborhood $V$ of $\left( c,b\right) $
such that 
\begin{equation}
d\left( \widetilde{x},\mathcal{S}\left( c,b\right) \right) \leq \kappa
d\left( \left( \widetilde{c},\widetilde{b}\right) ,\left( c,b\right) \right) 
\text{ for all }\left( \widetilde{c},\widetilde{b}\right) \in V\text{ and
all }\widetilde{x}\in \mathcal{S}\left( \widetilde{c},\widetilde{b}\right) .
\label{lipusc def}
\end{equation}%
This property will be used to characterize the Aubin property of $\mathcal{S}
$ and determine its Lipschitz modulus in this section. To do this we need to
introduce the Hausdorff lower semicontinuity property: $\mathcal{S}$ is said
to be Hausdorff lower semicontinuous at $\left( c,b\right) \in \limfunc{dom}%
\mathcal{S}$ if 
\begin{equation*}
\lim_{\left( \widetilde{c},\widetilde{b}\right) \rightarrow \left(
c,b\right) }e\left( \mathcal{S}\left( c,b\right) ,\mathcal{S}\left( 
\widetilde{c},\widetilde{b}\right) \right) =0,
\end{equation*}%
where $e\left( \mathcal{S}\left( c,b\right) ,\mathcal{S}\left( \widetilde{c},%
\widetilde{b}\right) \right) :=\sup_{x\in \mathcal{S}\left( c,b\right)
}d\left( x,\mathcal{S}\left( \widetilde{c},\widetilde{b}\right) \right) $ is
the Hausdorff excess of $\mathcal{S}\left( c,b\right) $ over $\mathcal{S}%
\left( \widetilde{c},\widetilde{b}\right) .$ In these terms, (\ref{lipusc
def}) reads as $e\left( \mathcal{S}\left( \widetilde{c},\widetilde{b}\right)
,\mathcal{S}\left( c,b\right) \right) \leq \kappa d\left( \left( \widetilde{c%
},\widetilde{b}\right) ,\left( c,b\right) \right) $ for all $\left( 
\widetilde{c},\widetilde{b}\right) \in V.$ At this moment we need some
additional technical results.

\begin{lem}
For any $x\in \mathbb{R}^{n}$ we have $x^{\prime }Qx=0\Leftrightarrow
Qx=0_{n}.$
\end{lem}

\begin{dem}
The proof is similar to the beginning of that of Proposition \ref{Prop M_D}.
Since $Q$ is assumed to be symmetric and positive semidefinite, it has a
(not necessarily unique) Cholesky decomposition $LL^{\prime }$. Hence, 
\begin{equation*}
0=x^{\prime }Qx=x^{\prime }LL^{\prime }x=\left( L^{\prime }x\right) ^{\prime
}L^{\prime }x\Rightarrow L^{\prime }x=0_{n}\Rightarrow Qx=0_{n}\Rightarrow
x^{\prime }Qx=0.
\end{equation*}
\end{dem}

If $\vartheta \left( c,b\right) $ denotes the optimal value of problem%
\begin{equation*}
P\left( c,b\right) :\text{minimize }\frac{1}{2}x^{\prime }Qx+c^{\prime }x%
\text{ such that }Ax\leq b,
\end{equation*}%
then, for any $\alpha \geq \vartheta \left( c,b\right) $ we can consider the 
$\alpha $-(sub)level set 
\begin{equation*}
\mathcal{L}\left( c,b,\alpha \right) :=\left\{ x\in \mathbb{R}^{n}:Ax\leq b,~%
\frac{1}{2}x^{\prime }Qx+c^{\prime }x\leq \alpha \right\} ,
\end{equation*}%
which is a nonempty closed convex set. The following lemma shows that the
recession cone of $\mathcal{L}\left( c,b,\alpha \right) $ does not depend on 
$\left( b,\alpha \right) .$ Recall from \cite[p. 61]{Rock70} that the
recession cone of a convex set $C,$ denoted as $0^{+}C$, is the set of all $%
u\in \mathbb{R}^{n}$ such that $z+\lambda u\in C$ for all $z\in C$ and all $%
\lambda \geq 0.$

\begin{lem}
\label{Lem recession indep RHS}Let $\left( c,b\right) \in \limfunc{dom}%
\mathcal{S}$ and $\alpha \geq \vartheta \left( c,b\right) $. Then 
\begin{equation*}
0^{+}\mathcal{L}\left( c,b,\alpha \right) =K:=\left\{ u\in \mathbb{R}%
^{n}:Au\leq 0_{m},~c^{\prime }u\leq 0,~Qu=0_{n}\right\} .
\end{equation*}
\end{lem}

\begin{dem}
It is clear that $K\subset 0^{+}\mathcal{L}\left( c,b,\alpha \right) .$ To
see the converse inclusion take any $u\in 0^{+}\mathcal{L}\left( c,b,\alpha
\right) .$ The relations $Au\leq 0_{m}$ and $c^{\prime }u\leq 0$ are direct
(see \cite[p. 62]{Rock70}). If $z\in \mathcal{L}\left( c,b,\alpha \right) $
and $\lambda \geq 0$ we have 
\begin{eqnarray*}
2\alpha &\geq &\left( z+\lambda u\right) ^{\prime }Q\left( z+\lambda
u\right) +2c^{\prime }\left( z+\lambda u\right) \\
&=&z^{\prime }Qz+2c^{\prime }z+2\lambda \left( z^{\prime }Qu+c^{\prime
}u\right) +\lambda ^{2}u^{\prime }Qu,
\end{eqnarray*}%
and dividing both members by $\lambda ^{2}$ and letting $\lambda \rightarrow
+\infty $ we get $u^{\prime }Qu\leq 0.$ Since $Q$ is positive semidefinite
we have $u^{\prime }Qu=0,$ which is equivalent to $Qu=0_{n}$ because of the
previous lemma.
\end{dem}

The following result comes from standard arguments in parametric
optimization (see, e. g., \cite[Section 5.5]{BGKKT82}). For the sake of
completeness we include a proof here.

\begin{prop}
\label{Prop S Hlsc}Let $\left( c,b\right) \in \mathbb{R}^{n}\times \mathbb{R}%
^{m}$ be such that $\mathcal{S}\left( c,b\right) $ is a singleton and SCQ
holds at $b$. Then $\mathcal{S}$ is Hausdorff lower semicontinuous at $%
\left( c,b\right) .$
\end{prop}

\begin{dem}
Write $\mathcal{S}\left( c,b\right) =\{x\}$ and take any $\alpha >\vartheta
\left( c,b\right) =\frac{1}{2}x^{\prime }Qx+c^{\prime }x.$ Let $y$ be a
Slater element for $b,$ that is, $a_{i}^{\prime }y<b_{i}$ for $i=1,...,m.$
Then it is clear that $y^{\mu }:=\left( 1-\mu \right) x+\mu y$ is also a
Slater element for $b$ whenever $0<\mu \leq 1,$ and if such a $\mu $ is
small enough we also have $\frac{1}{2}\left( y^{\mu }\right) ^{\prime
}Qy^{\mu }+c^{\prime }y^{\mu }<\alpha .$ Fix $\mu $ satisfying these
conditions. Then there exists $\delta >0$ such that $\left\Vert \widetilde{b}%
-b\right\Vert _{\infty }<\delta $ and $\left\Vert \widetilde{c}-c\right\Vert
_{\ast }<\delta $ imply 
\begin{equation*}
Ay^{\mu }\leq \widetilde{b}\text{ and }\frac{1}{2}\left( y^{\mu }\right)
^{\prime }Qy^{\mu }+\widetilde{c}^{\prime }y^{\mu }<\alpha .
\end{equation*}%
This yields $y^{\mu }\in \mathcal{L}\left( \widetilde{c},\widetilde{b}%
,\alpha \right) ,$ and hence this set is nonempty.

On the other hand, The boundedness of $\mathcal{S}\left( c,b\right) $
entails that $0^{+}\mathcal{S}\left( c,b\right) =\left\{ 0_{n}\right\} $
(see \cite[Theorem 8.4]{Rock70}) and then $0^{+}\mathcal{L}\left( c,b,\alpha
\right) =\left\{ 0_{n}\right\} $ according to Lemma \ref{Lem recession indep
RHS} and the preceding comments to that Lemma. Let us see that $\delta $
above can be chosen small enough to guarantee $0^{+}\mathcal{L}\left( 
\widetilde{c},\widetilde{b},\alpha \right) =\left\{ 0_{n}\right\} $ whenever 
$\left\Vert \widetilde{b}-b\right\Vert _{\infty }<\delta $ and $\left\Vert 
\widetilde{c}-c\right\Vert _{\ast }<\delta .$ Assume to the contrary that we
can take a sequence $\left\{ \left( c^{r},b^{r}\right) \right\} _{r\in 
\mathbb{N}}$ converging to $\left( c,b\right) $ and for each $r$ a certain $%
u^{r}\in 0^{+}\mathcal{L}\left( c^{r},b^{r},\alpha \right) $ with $%
\left\Vert u^{r}\right\Vert =1$. Without loss of generality we may assume
that $\left\{ u^{r}\right\} _{r\in \mathbb{N}}$ converges to certain $u\in 
\mathbb{R}^{n}$ with $\left\Vert u\right\Vert =1$. Then, the description of
the recession cones of the level sets given in Lemma \ref{Lem recession
indep RHS} clearly entails $u\in 0^{+}\mathcal{L}\left( c,b,\alpha \right) ,$
a contradiction.

Applying again \cite[Theorem 8.4]{Rock70} we conclude the boundedness of $%
\mathcal{L}\left( \widetilde{c},\widetilde{b},\alpha \right) ,$ and
therefore (applying Weierstrass Theorem) the nonemptiness of $\mathcal{S}%
\left( \widetilde{c},\widetilde{b}\right) $ whenever $\left\Vert \widetilde{b%
}-b\right\Vert _{\infty }<\delta $ and $\left\Vert \widetilde{c}%
-c\right\Vert _{\ast }<\delta .$ Now, fixing arbitrarily $\varepsilon >0,$
the Lipschitz upper semicontinuity of $\mathcal{S}$ at $\left( c,b\right) $
(see (\ref{lipusc def})) ensures that, if $\delta $ is taken small enough, $%
\left\Vert \widetilde{b}-b\right\Vert _{\infty }<\delta $ and $\left\Vert 
\widetilde{c}-c\right\Vert _{\ast }<\delta $ entail 
\begin{equation}
\emptyset \neq \mathcal{S}\left( \widetilde{c},\widetilde{b}\right) \subset
x+\varepsilon \mathbb{B},  \label{8bis}
\end{equation}%
where $\mathbb{B}$ denotes the closed unit ball with respect to norm $%
\left\Vert \cdot \right\Vert ,$ which entails 
\begin{equation*}
e\left( \mathcal{S}\left( c,b\right) ,\mathcal{S}\left( \widetilde{c},%
\widetilde{b}\right) \right) =d\left( x,\mathcal{S}\left( \widetilde{c},%
\widetilde{b}\right) \right) \leq \varepsilon .
\end{equation*}
\end{dem}

Next we establish the main result of the paper under our standing assumption
that SCQ holds at $\overline{b}.$ This result characterizes the Aubin
property of $\mathcal{S}$ at a point $\left( \left( \overline{c},\overline{b}%
\right) ,\overline{x}\right) $ of its graph and provides the Lipschitz
modulus at this point. According to Proposition \ref{Prop1}, the hypothesis $%
\mathcal{S}\left( \overline{c},\overline{b}\right) =\left\{ \overline{x}%
\right\} $ is not restrictive at all.

\begin{theo}
\label{Th main}Assume $\mathcal{S}\left( \overline{c},\overline{b}\right)
=\left\{ \overline{x}\right\} .$ The following conditions are equivalent:

$\left( i\right) $ $\mathcal{S}$ has the Aubin property at $\left( \left( 
\overline{c},\overline{b}\right) ,\overline{x}\right) ;$

$\left( ii\right) $ $M_{D}$ is nonsingular for all $D\in \mathcal{M}_{\bar{c}%
,\bar{b}}\left( \bar{x}\right) ;$

$\left( iii\right) $ $M_{D}$ is nonsingular for all $D\in \mathcal{L}_{\bar{c%
},\bar{b}}\left( \bar{x}\right) .$

Moreover, under these equivalent conditions we have 
\begin{equation}
\mathrm{lip}\mathcal{S}\left( \left( \overline{c},\overline{b}\right) ,%
\overline{x}\right) =\max_{D\in \mathcal{L}_{\bar{c},\bar{b}}\left( \bar{x}%
\right) }\left\Vert \left( 
\begin{array}{cc}
I_{n} & 0_{n\times \left\vert D\right\vert }%
\end{array}%
\right) M_{D}^{-1}\right\Vert .  \label{6bis}
\end{equation}
\end{theo}

\begin{dem}
The implication $\left( iii\right) \Rightarrow \left( ii\right) $ is evident
from the fact that $\mathcal{M}_{\bar{c},\bar{b}}\left( \bar{x}\right)
\subset \mathcal{L}_{\bar{c},\bar{b}}\left( \bar{x}\right) ,$ whereas the
converse implication comes from Corollary \ref{Cor_D0}, taking into account
that every element of $\mathcal{L}_{\bar{c},\bar{b}}\left( \bar{x}\right) $
contains an element of $\mathcal{M}_{\bar{c},\bar{b}}\left( \bar{x}\right) $.

$\left( i\right) \Rightarrow \left( ii\right) $ follows directly from
Corollary \ref{CorS pL implies D nonsingular}.

$\left( iii\right) \Rightarrow \left( i\right) .$ From Proposition \ref{Prop
S Hlsc} it is clear that $\mathcal{S}\left( c,b\right) $ is nonempty for $%
\left( c,b\right) $ close enough to $\left( \overline{c},\overline{b}\right)
.$ Let us see that $\left( iii\right) $ implies the existence of a
neighborhood $V$ of $\left( \overline{c},\overline{b}\right) $ such that $%
\mathcal{S}$ is single-valued on $V.$ Assume to the contrary that there
exists a sequence $\left\{ \left( c^{r},b^{r}\right) \right\} _{r\in \mathbb{%
N}}$ converging to $\left( \overline{c},\overline{b}\right) $ such that, for
all $r,$ $\mathcal{S}\left( c^{r},b^{r}\right) $ is not a singleton. For
each $r$ pick $x^{r}\in \mathcal{S}\left( c^{r},b^{r}\right) $ and $D_{r}\in 
\mathcal{M}_{c^{r},b^{r}}\left( x^{r}\right) $ and assume (by taking a
subsequence, if necessary, and recalling that $I$ is finite) that the
sequence $\left\{ D_{r}\right\} _{r\in \mathbb{N}}$ is constant, say $%
D_{r}=D $ for each $r.$ From Proposition \ref{Prop Limsup minimimals}$\left(
i\right) $ it is clear that $D\in \mathcal{L}_{\bar{c},\bar{b}}\left( \bar{x}%
\right) $ and that $x^{r}\in \mathcal{S}_{D}\left( c^{r},b_{D}^{r}\right) $,
because of the KKT conditions given by $D\in \mathcal{M}_{c^{r},b^{r}}\left(
x^{r}\right) .$ By $\left( iii\right) $ we have that $M_{D}$ is nonsingular
and, then, equivalence $\left( iv\right) \Leftrightarrow \left( v\right) $
in Theorem \ref{Theorem_SD}, applied to $\left( \left( c^{r},b^{r}\right)
,x^{r}\right) \in \limfunc{gph}\mathcal{S},$ gives $\mathcal{S}_{D}\left(
c^{r},b_{D}^{r}\right) =\{x^{r}\}$ and, recalling Remark \ref{Rem_S_SD}$%
\left( iii\right) $, 
\begin{equation*}
\mathcal{S}\left( c^{r},b^{r}\right) =\mathcal{S}_{D}\left(
c^{r},b_{D}^{r}\right) \cap \mathcal{F}\left( b^{r}\right) =\{x^{r}\},
\end{equation*}%
contradicting our current assumption. Hence, we have proved that $\mathcal{S}
$ is single-valued in some neighborhood of $\left( \overline{c},\overline{b}%
\right) .$

Under the assumption of $\left( iii\right) ,$ the aimed condition $\left(
i\right) $ follows immediately if we prove inequality `$\leq $' in (\ref%
{6bis}). To do this, take 
\begin{equation*}
\kappa :=\max_{D\in \mathcal{L}_{\bar{c},\bar{b}}\left( \bar{x}\right)
}\left\Vert \left( 
\begin{array}{cc}
I_{n} & 0_{n\times \left\vert D\right\vert }%
\end{array}%
\right) M_{D}^{-1}\right\Vert ,
\end{equation*}%
and let us find a neighborhood $V\times U$ of $\left( \left( \overline{c},%
\overline{b}\right) ,\overline{x}\right) $ where $\kappa $ is a Lipschitz
constant for $\mathcal{S},$ that is, (\ref{eq_def_Lips}) holds in $V\times
U. $ To do this we will apply \cite[Theorem 2.1]{Li94}, which guarantees
that if a set-valued mapping from a convex subset $C\subset \mathbb{R}^{k}$
to the subsets of $\mathbb{R}^{n}$ is Hausdorff lower semicontinuous on $C$
and Lipschitz upper semicontinuous on $C$ with the same Lipschitz upper
semicontinuity constant $\kappa $ at all the points of $C,$ then $\kappa $
becomes a Lipschitz constant on $C.$

Take an open and convex neighborhood $V\times U$ of $\left( \left( \overline{%
c},\overline{b}\right) ,\overline{x}\right) $ such that:

$\cdot $ the SCQ holds at all $b\in \mathbb{R}^{m}$ with $\left( c,b\right)
\in V$ for some $c\in \mathbb{R}^{n},$

$\cdot $ $\mathcal{S}$ is single-valued on $V$, say $\mathcal{S}\left(
c,b\right) =\left\{ x\left( c,b\right) \right\} $ for $\left( c,b\right) \in
V,$

$\cdot $ $\mathcal{S}$ is Hausdorff lower semicontinuous on $V$ (apply
Proposition \ref{Prop S Hlsc}),

$\cdot $ $\mathcal{S}\left( V\right) :=\bigcup\nolimits_{\left( c,b\right)
\in V}\mathcal{S}\left( c,b\right) \subset U$ (possible because of the
previous item),

$\cdot $ (\ref{eq_MinL}) holds for $V\times U$.

Fix any $\left( c^{0},b^{0}\right) \in V$ and, applying Propositions \ref%
{Prop Limsup minimimals}$\left( i\right) $ and \ref{Prop S Hlsc} at $(\left(
c^{0},b^{0}\right) ,\allowbreak x\left( c^{0},b^{0}\right) ),$ take a
neighborhood $V_{0}\times U_{0}\subset V\times U$ of $\left( \left(
c^{0},b^{0}\right) ,x\left( c^{0},b^{0}\right) \right) $ such that $x\left(
V_{0}\right) :=\left\{ x\left( c,b\right) ,~\left( c,b\right) \in
V_{0}\right\} \subset U_{0}$ and 
\begin{equation*}
\mathcal{M}_{c,b}\left( x\left( c,b\right) \right) \subset \mathcal{L}%
_{c^{0},b^{0}}\left( x\left( c^{0},b^{0}\right) \right) \text{ for all }%
\left( c,b\right) \in V_{0}.
\end{equation*}%
For any $\left( c,b\right) \in V_{0}$ pick $D\in \mathcal{M}_{c,b}\left(
x\left( c,b\right) \right) .$ Proposition \ref{Prop Limsup minimimals}$%
\left( i\right) ,\left( iii\right) $ yields 
\begin{equation*}
D\in \mathcal{L}_{c^{0},b^{0}}\left( x\left( c^{0},b^{0}\right) \right)
\subset \mathcal{L}_{\bar{c},\bar{b}}\left( \bar{x}\right) .
\end{equation*}%
The KKT optimality conditions yield the existence of $\lambda ^{0},\lambda
\in \mathbb{R}_{+}^{D}$ such that 
\begin{equation*}
M_{D}\dbinom{x\left( c^{0},b^{0}\right) }{\lambda ^{0}}=\dbinom{-c^{0}}{%
b_{D}^{0}}\text{ and }M_{D}\dbinom{x\left( c,b\right) }{\lambda }=\dbinom{-c%
}{b_{D}}.
\end{equation*}%
Then, the nonsingularity of $M_{D}$ yields 
\begin{equation*}
x\left( c^{0},b^{0}\right) =\left( 
\begin{array}{cc}
I_{n} & 0_{n\times \left\vert D\right\vert }%
\end{array}%
\right) M_{D}^{-1}\dbinom{-c^{0}}{b_{D}^{0}}\text{ and }x\left( c,b\right)
=\left( 
\begin{array}{cc}
I_{n} & 0_{n\times \left\vert D\right\vert }%
\end{array}%
\right) M_{D}^{-1}\dbinom{-c}{b_{D}}.
\end{equation*}%
Accordingly, from the definition of matrix norm,%
\begin{eqnarray*}
\left\Vert x\left( c,b\right) -x\left( c^{0},b^{0}\right) \right\Vert &\leq
&\left\Vert \left( 
\begin{array}{cc}
I_{n} & 0_{n\times \left\vert D\right\vert }%
\end{array}%
\right) M_{D}^{-1}\right\Vert \left\Vert \dbinom{-c}{b_{D}}-\dbinom{-c^{0}}{%
b_{D}^{0}}\right\Vert \\
&\leq &\kappa \left\Vert \dbinom{c}{b_{D}}-\dbinom{c^{0}}{b_{D}^{0}}%
\right\Vert \leq \kappa \left\Vert \dbinom{c}{b}-\dbinom{c^{0}}{b^{0}}%
\right\Vert ,
\end{eqnarray*}%
taking into account the symmetry of the norm, $\max \left\{ \left\Vert \cdot
\right\Vert _{\ast },\left\Vert \cdot \right\Vert _{\infty }\right\} ,$
being used in $\mathbb{R}^{n}\times \mathbb{R}^{D}.$ Therefore, $\kappa $ is
a Lipschitz upper semicontinuity constant of $\mathcal{S}$ at any $\left(
c^{0},b^{0}\right) \in V.$

Finally, we can apply in $V$ the variant of \cite[Theorem 2.1]{Li94}
referred to above to conclude (\ref{eq_def_Lips}), namely 
\begin{equation*}
\left\Vert x\left( c^{1},b^{1}\right) -x\left( c^{2},b^{2}\right)
\right\Vert \leq \kappa \left\Vert \dbinom{c^{1}}{b^{1}}-\dbinom{c^{2}}{b^{2}%
}\right\Vert
\end{equation*}%
for all $\left( c^{1},b^{1}\right) ,\left( c^{2},b^{2}\right) \in V.$ This
entails $\mathrm{lip}\mathcal{S}\left( \left( \overline{c},\overline{b}%
\right) ,\overline{x}\right) \leq \kappa $, i.e., inequality `$\leq $' in (%
\ref{6bis}) holds; hence $\left( i\right) $ holds.

Let us show that the converse inequality also holds. From the definition of
Lipschitz modulus we immediately deduce 
\begin{equation}
\mathrm{lip}\mathcal{S}\left( \left( \overline{c},\overline{b}\right) ,%
\overline{x}\right) \geq \underset{_{\substack{ \left( \left( c,b\right)
,x\right) \rightarrow \left( \left( \overline{c},\overline{b}\right) ,%
\overline{x}\right)  \\ \left( \left( c,b\right) ,x\right) \in \mathrm{gph}%
\mathcal{S}}}}{\lim \sup }\mathrm{lip}\mathcal{S}\left( \left( c,b\right)
,x\right) .  \label{7bis}
\end{equation}%
More in detail, take any $\kappa >\mathrm{lip}\mathcal{S}\left( \left( 
\overline{c},\overline{b}\right) ,\overline{x}\right) ,$ the latter being
finite under $\left( i\right) $. Then there exists an open neighborhood $%
V\times U$ of $\left( \left( \overline{c},\overline{b}\right) ,\overline{x}%
\right) $ such that (\ref{eq_def_Lips}) holds whenever $x^{2}\in \mathcal{S}%
\left( c^{2},b^{2}\right) \cap U$ and $\left( c^{1},b^{1}\right) ,\left(
c^{2},b^{2}\right) \in V.$ Take any $\left( \left( c,b\right) ,x\right) \in
\left( V\times U\right) \cap \mathrm{gph}\mathcal{S}$ and any neighborhood $%
V_{0}\times U_{0}$ of $\left( \left( c,b\right) ,x\right) $ contained in $%
V\times U.$ Then, the particularization of (\ref{eq_def_Lips}) to $%
V_{0}\times U_{0}$ ensures that $\kappa \geq \mathrm{lip}\mathcal{S}\left(
\left( c,b\right) ,x\right) $. This entails (\ref{7bis}).

Now, given any $D\in \mathcal{L}_{\bar{c},\bar{b}}\left( \bar{x}\right) ,$
lets us apply Proposition \ref{Prop Limsup minimimals}$\left( ii\right) $ to
write $\left( \left( \overline{c},\overline{b}\right) ,\overline{x}\right)
=\lim\limits_{r\rightarrow \infty }\left( \left( c^{r},b^{r}\right)
,x^{r}\right) $ with $\left( \left( c^{r},b^{r}\right) ,x^{r}\right) \in 
\mathrm{gph}\mathcal{S}$ for all $r,$ and $D\in \mathcal{M}%
_{c^{r},b^{r}}\left( x^{r}\right) .$ Then (\ref{7bis}) together with
Proposition \ref{Prop>=} and Corollary \ref{Cor_lipSD} entail 
\begin{eqnarray*}
\mathrm{lip}\mathcal{S}\left( \left( \overline{c},\overline{b}\right) ,%
\overline{x}\right) &\geq &\underset{r\rightarrow \infty }{\lim \sup }\,%
\mathrm{lip}\mathcal{S}\left( \left( c^{r},b^{r}\right) ,x^{r}\right) \geq 
\underset{r\rightarrow \infty }{\lim \sup }\,\mathrm{lip}\mathcal{S}%
_{D}\left( \left( c^{r},b_{D}^{r}\right) ,x^{r}\right) \\
&=&\left\Vert \left( 
\begin{array}{cc}
I_{n} & 0_{n\times \left\vert D\right\vert }%
\end{array}%
\right) M_{D}^{-1}\right\Vert .
\end{eqnarray*}%
Observe that the last $\lim \sup $ affects a constant sequence. Hence, we
have proved inequality `$\geq $' in (\ref{6bis}).
\end{dem}

\begin{cor}
Assume $\mathcal{S}\left( \overline{c},\overline{b}\right) =\left\{ 
\overline{x}\right\} $ and let $D_{0}\in \mathcal{L}_{\bar{c},\bar{b}}\left( 
\bar{x}\right) .$ The following conditions are equivalent:

$\left( i\right) $ $\mathcal{S}_{D_{0}}$ has the Aubin property at $\left(
\left( \overline{c},\overline{b}_{D_{0}}\right) ,\overline{x}\right) ;$

$\left( ii\right) $ $\mathcal{S}_{D}$ has the Aubin property at $\left(
\left( \overline{c},\overline{b}_{D}\right) ,\overline{x}\right) $ for all $%
D\in \mathcal{M}_{\bar{c},\bar{b}}\left( \bar{x}\right) $ with $D\subset
D_{0};$

$\left( iii\right) $ $\mathcal{S}_{D}$ has the Aubin property at $\left(
\left( \overline{c},\overline{b}_{D}\right) ,\overline{x}\right) $ for all $%
D\in \mathcal{L}_{\bar{c},\bar{b}}\left( \bar{x}\right) $ with $D\subset
D_{0};$

$\left( iv\right) $ $M_{D}$ is nonsingular for all $D\in \mathcal{L}_{\bar{c}%
,\bar{b}}\left( \bar{x}\right) $ with $D\subset D_{0}.$

Moreover, under these equivalent conditions we have 
\begin{equation}
\mathrm{lip}\mathcal{S}_{D_{0}}\left( \left( \overline{c},\overline{b}%
_{D_{0}}\right) ,\overline{x}\right) =\max_{D\in \mathcal{L}_{\bar{c},\bar{b}%
}\left( \bar{x}\right) \cap 2^{D_{0}}}\left\Vert \left( 
\begin{array}{cc}
I_{n} & 0_{n\times \left\vert D\right\vert }%
\end{array}%
\right) M_{D}^{-1}\right\Vert .  \label{eq11}
\end{equation}
\end{cor}

\begin{dem}
The result comes from applying the previous theorem to problem $%
P_{D_{0}}\left( \overline{c},\overline{b}_{D_{0}}\right) $ playing the role
of the whole $P\left( \overline{c},\overline{b}\right) .$ More specifically,
condition $\left( ii\right) $ in the present corollary is the same as
condition $\left( ii\right) $ in Theorem \ref{Th main} as a consequence of $%
\left( i\right) \Leftrightarrow \left( v\right) $ in Theorem \ref{Theorem_SD}%
. Hence $\left( i\right) \Leftrightarrow \left( ii\right) $ in the present
corollary holds. Implication $\left( iii\right) \Rightarrow \left( ii\right) 
$ is obvious, whereas $\left( ii\right) \Rightarrow \left( iii\right) $
comes to apply $\left( ii\right) \Rightarrow \left( i\right) $ to each
problem $P_{D}\left( \overline{c},\overline{b}_{D}\right) $ with $D\in 
\mathcal{L}_{\bar{c},\bar{b}}\left( \bar{x}\right) $ and $D\subset D_{0}.$
In addition, $\left( ii\right) \Leftrightarrow \left( iv\right) $ comes from
Corollary \ref{Cor_D0}. Finally, the last assertion of the corollary comes
from (\ref{6bis}) applied to problem $P_{D_{0}}\left( \overline{c},\overline{%
b}_{D_{0}}\right) .$
\end{dem}

The following example shows that the maximum in (\ref{eq11}) may not be
attained at $D_{0}.$

\begin{exa}
\label{Exa2}\emph{Consider the parameterized problem in} $\mathbb{R}^{2},$ 
\emph{endowed with the Euclidean norm, denoted as }$\left\Vert \cdot
\right\Vert _{2},$%
\begin{equation*}
\begin{tabular}{lll}
$P\left( c,b\right) :$ & \textrm{Minimize} & $\frac{1}{2}x^{\prime
}Qx+c^{\prime }x$ \\ 
& \textrm{subject to} & $Ax\leq b,$%
\end{tabular}%
\end{equation*}%
\emph{where }$Q=\left( 
\begin{array}{cc}
0 & 0 \\ 
0 & \alpha%
\end{array}%
\right) $\emph{\ with }$\alpha >0,$ $A=\left( 
\begin{array}{cc}
-1 & 1 \\ 
-1 & -1 \\ 
-1 & 0%
\end{array}%
\right) ,$\emph{\ nominal parameters }$\overline{c}=\binom{1}{0}$ \emph{and }%
$\overline{b}=0_{3},$ \emph{and nominal solution }$\overline{x}=0_{2}.$
\end{exa}

\noindent Then we have 
\begin{equation*}
\mathcal{L}_{\overline{c},\overline{b}}\left( \overline{x}\right) =\left\{
\left\{ 1,2\right\} ,\left\{ 3\right\} ,\left\{ 1,3\right\} ,\left\{
2,3\right\} \right\} ,
\end{equation*}%
whereas $\mathcal{M}_{\bar{c},\bar{b}}\left( \bar{x}\right) =\left\{ \left\{
1,2\right\} ,\left\{ 3\right\} \right\} .$ We also have the following table:%
\begin{equation*}
\begin{tabular}{||c||c|c|c|c||}
\hline\hline
$D$ & $\left\{ 1,2\right\} $ & $\left\{ 3\right\} $ & $\left\{ 1,3\right\} $
& $\left\{ 2,3\right\} $ \\ \hline
$\left\Vert \left( 
\begin{array}{cc}
I_{n} & 0_{n\times \left\vert D\right\vert }%
\end{array}%
\right) M_{D}^{-1}\right\Vert $ & $1$ & $\sqrt{1+\left( 1/\alpha \right) ^{2}%
}$ & $\sqrt{5}$ & $\sqrt{5}$ \\ \hline\hline
\end{tabular}%
\end{equation*}

According to Theorem \ref{Th main}, we have 
\begin{equation*}
\mathrm{lip}\mathcal{S}\left( \left( \overline{c},\overline{b}\right) ,%
\overline{x}\right) =\sqrt{1+\left( \max \left\{ 2,1/\alpha \right\}
^{2}\right) }.
\end{equation*}%
In particular, for $\alpha =1/3$ and $D_{0}=\left\{ 1,3\right\} $ the
maximum in (\ref{eq11}) is not attained at $D_{0},$ but at the proper subset 
$\left\{ 3\right\} .$

\section{Conclusions}

In this section we summarize the main results of the paper and come back to
the computation of the Lipschitz modulus of the metric projection (\ref{eq_P}%
) as its initial motivation. The metric projection corresponds to the
particular case of (\ref{eq_S}) given by $Q=I_{n}$ (the identity matrix in $%
\mathbb{R}^{n\times n}$) and $c=-z,$ with $z$ being the point to be
projected on the convex polyhedron (\ref{eq_F}). In this paper we analyze
the convex quadratic problem when $Q$ is any given positive semidefinite
matrix. This includes the linear programming setting, widely studied in the
literature, as the particular case $Q=0_{n\times n}.$

The convex quadratic setting presents notable differences with respect to
the linear one. Clearly, the optimal value of the former, when finite, needs
not be attained at an extreme point (if any) of the feasible set. While in
the linear setting the Aubin property of the argmin mapping $\mathcal{S}$ is
characterized by the N\"{u}rnberger condition (see Section 2.2), in the
convex quadratic setting this condition is only sufficient for the Aubin
property of $\mathcal{S},$ which may be characterized in terms of the
minimal sets of KKT indices; see Theorem \ref{Th main}$\left( ii\right) .$
However, these minimal sets are not enough to compute the Lipschitz modulus
of $\mathcal{S};$ see Proposition \ref{Prop>=} and Example \ref{Exa1}. In
order to obtain a point-based formula (depending only on the nominal point
and parameters, not involving elements in a neighborhood) we need to
consider the so-called extended family of KKT subsets of indices; see (\ref%
{eq_M in L}). The mentioned Lipschitz modulus, see (\ref{6bis}), is also
given in Theorem \ref{Th main} (the main result of the paper).

Coming back to the metric projection, the fact that $Q=I_{n}$ is nonsingular
entails, by taking \thinspace $D_{0}=\emptyset $ in Corollary \ref{Cor_D0}
and applying Theorem \ref{Th main}, that the metric projection function $%
\mathcal{P}$ in (\ref{eq_P}) has the Aubin property at any $\left( \overline{%
z},\overline{b}\right) \in \mathbb{R}^{n}\times \mathbb{R}^{m}$ such that
the Slater condition holds at $\overline{b};$ see our standing assumption at
the end of Section 2.1 and recall the bibliographic comments at the
beginning of Section 1. The present paper allows us to go further and
provide the following point-based expression for the Lipschitz modulus of
the metric projection (the vertical bar is to emphasize that we are
considering partitioned matrices):%
\begin{equation}
\mathrm{lip}\mathcal{P}\left( \left( \overline{z},\overline{b}\right) ,%
\overline{x}\right) =\max_{D\in \mathcal{L}_{\bar{z},\bar{b}}\left( \bar{x}%
\right) }\left\Vert \left( I_{n}-A_{D}^{\prime }\left( A_{D}A_{D}^{\prime
}\right) ^{-1}A_{D}\mid A_{D}^{\prime }\left( A_{D}A_{D}^{\prime }\right)
^{-1}\right) \right\Vert ,  \label{lip proj}
\end{equation}%
where $\overline{x}$ is the (unique) projection of $\overline{z}$ on $%
\mathcal{F}\left( \overline{b}\right) $ and 
\begin{equation*}
\mathcal{L}_{\bar{z},\bar{b}}\left( \bar{x}\right) =\left\{ D\subset I_{%
\overline{b}}\left( \overline{x}\right) \left\vert 
\begin{array}{l}
\left\{ a_{i},i\in D\right\} \text{ is linearly independent} \\ 
\text{and }\overline{z}-\overline{x}\in \mathrm{cone}\left\{ a_{i},i\in
D\right\} 
\end{array}%
\right. \right\} .
\end{equation*}%
Here (\ref{lip proj}) comes directly from (\ref{6bis}) and the fact that,
for $Q=I_{n},$ the inverse matrix $M_{D}^{-1}$ can be easily computed,
taking into account that $A_{D}A_{D}^{\prime }$ is nonsingular because $A_{D}
$ has full-row-rank. Observe that we have adapted the definition of $%
\mathcal{L}_{\bar{c},\bar{b}}\left( \bar{x}\right) $ to the case $\overline{c%
}=-\overline{z}.\medskip $

\textbf{Declaration: }There are no competing interests.


\begin{thebibliography}{99}
\bibitem{Ab79} T. ABATZOGLOU, \emph{The Lipschitz continuity of the metric
projection}. J. Approx. Theory 26, 212-218 (1979)

\bibitem{BaGo12} M. V. BALASHOV, M. O. GOLUBEV, \emph{About the Lipschitz
property of the metric projection in the Hilbert space}. J. Math. Anal.
Appl. 394, 545-551 (2012)

\bibitem{BGKKT82} B. BANK, J. GUDDAT, D. KLATTE, B. KUMMER, K. TAMMER, \emph{%
Non-Linear Parametric Optimization}, Akademie-Verlag, Berlin, 1982, and Birkh%
\"{a}user, Basel, 1983.

\bibitem{BeRu19} E. M. BERNARCZUK, K. E. RUTKOWSKI, \emph{On Lipschitz
continuity of projections onto polyhedral moving sets}.
arXiv:1909.13715v2[math.OC], 2019.

\bibitem{CDLP05} M. J. C\'{A}NOVAS, A.L. DONTCHEV, M. A. L\'{O}PEZ, J.
PARRA, \emph{Metric regularity of semi-infinite constraint systems}. Math.
Program. B 104, 329--346 (2005)

\bibitem{CGP08} M. J. C\'{A}NOVAS, F. J. G\'{O}MEZ-SENENT, J. PARRA,\emph{\
On the Lipschitz modulus of the argmin mapping in linear semi-infinite
optimization. }Set-Valued Anal. Set Valued Anal. 16, 511--538 (2008)

\bibitem{CHLP08} M. J. C\'{A}NOVAS, A. HANTOUTE, M. A. L\'{O}PEZ, J. PARRA, 
\emph{Lipschitz modulus of the optimal set mapping in convex semi-infinite
optimization via minimal subproblems}. Pacific Journal of Optimization 4,
411-422 (2008)

\bibitem{chlp16} M. J. C\'{A}NOVAS, R. HENRION, M. A. L\'{O}PEZ, J. PARRA, 
\emph{Outer limit of subdifferentials and calmness moduli in linear and
nonlinear programming.} J. Optim. Theory Appl. 169, 925--952 (2016)

\bibitem{CKLP07} M. J. C\'{A}NOVAS, D. KLATTE, M. A. L\'{O}PEZ, J. PARRA, 
\emph{Metric regularity in convex semi-infinite optimization under canonical
perturbations}. SIAM J. Optim. 18, 717--732 (2007)

\bibitem{CLPT14} M. J. C\'{A}NOVAS, M. A. L\'{O}PEZ, J. PARRA, F. J. TOLEDO, 
\emph{Calmness of the feasible set mapping for linear inequality systems}.
Set-Valued Var. Anal. 22, 375--389 (2014)

\bibitem{DoRo} A. L. DONTCHEV, R. T. ROCKAFELLAR, \emph{Implicit Functions
and Solution Mappings: A View from Variational Analysis}, Springer, New
York, 2009.

\bibitem{libro} M. A. GOBERNA, M. A. L\'{O}PEZ, \emph{Linear Semi-Infinite
Optimization}, John Wiley \& Sons, Chichester (UK), 1998.

\bibitem{Ioffe17} A.D. IOFFE, \emph{Variational Analysis of Regular
Mappings, }Springer, Cham, 2017.

\bibitem{KYF04} C. KANZOW, N. YAMASHITA, M. FUKUSHIMA, \emph{%
Levenberg-Marquardt methods with strong local convergence properties f or
solving nonlinear equations with convex constraints}, J. Comput. Appl.
Math., 172, 375-397\emph{\ }(2004).

\bibitem{KlKu02} D. KLATTE, B. KUMMER, \emph{Nonsmooth Equations in
Optimization: Regularity, Calculus, Methods and Applications}, Nonconvex
Optim. Appl. 60, Kluwer Academic, Dordrecht, The Netherlands, 2002.

\bibitem{KlaKu05} D. KLATTE, B. KUMMER, \emph{Strong Lipschitz stability of
stationary solutions for nonlinear programs and variational inequalities}.%
\emph{\ }SIAM J. Optim. 16, 96-119 (2005).

\bibitem{KK09} D. KLATTE, B. KUMMER, \emph{Optimization methods and
stability of inclusions in Banach spaces}. Math. Program. Ser. B 117,
305--330 (2009)

\bibitem{LeeYen14} G. M. LEE, N. D. YEN, \emph{Coderivatives of
Karush-Kuhn-Tucker point set map and applications}. Nonlinear Anal. 95,
191-201 (2014)

\bibitem{Li94} W. LI, \emph{Sharp Lipschitz constants for basic optimal
solutions and basic feasible solutions of linear programs.} SIAM J. Control
Optim. 32, 140--153(1994)

\bibitem{Klatte25} D. KLATTE, \emph{Lipschitz stability for a class of
parametric optimization problems with polyhedral feasible set mapping.}
Optimization Online, https://optimization-online.org/?p=30863 (2025)

\bibitem{mor06a} B. S. MORDUKHOVICH,{\normalsize \ \emph{Variational
Analysis and Generalized Differentiation, I: Basic Theory}, }Springer,
Berlin, 2006{\normalsize . }

\bibitem{Robinson} S. M. ROBINSON, \emph{Some continuity properties of
polyhedral multifunctions}. Mathematical programming at Oberwolfach (Proc.
Conf., Math. Forschungsinstitut, Oberwolfach, 1979). Math. Programming Stud. 
\textbf{14}, 206--214 (1981)

\bibitem{Rock70} R. T. ROCKAFELLAR, {\normalsize \emph{Convex Analysis}, }%
Princeton University Press, Princeton, N.J., 1970.

\bibitem{rw} R. T. ROCKAFELLAR, R. J-B. WETS{\normalsize , \emph{Variational
Analysis}, }Springer, Berlin, 1998.

\bibitem{WYYZZ18} X. WANG, J.J. YE, X. YUAN, S. ZENG, J. ZHANG, \emph{%
Perturbation techniques for convergence analysis of proximal gradient method
and other first-order algorithms via variational analysis}. Ser-Valued Var.
Anal. \textbf{30}, 39-79 (2022)\P 

\bibitem{Yen95} N. D. YEN, \emph{Lipschitz continuity of solutions of
variational inequalities with a parametric polyhedral constraint}. Math.
Oper. Res. 20, 695-708 (1995)
\end{thebibliography}
\end{document}